\documentclass[preprint]{elsarticle}

\usepackage[english]{babel}
\usepackage{hyphenat}
\usepackage{acro}

\usepackage[a4paper, total={6.5in, 8.66in}]{geometry}
\usepackage{tocbibind}
\usepackage[export]{adjustbox}
\usepackage[framemethod=tikz]{mdframed}
\usepackage{float}

\usepackage{longtable}
\usepackage{booktabs}
\usepackage{multirow}
\usepackage{multicol}

\usepackage{amsmath}
\usepackage{amssymb}
\usepackage{mathtools}
\usepackage{amsthm}
\usepackage[per-mode=symbol,
            inter-unit-product=\cdot,
            separate-uncertainty=true,
            print-unity-mantissa=false,
            table-alignment-mode=none,
            table-alignment=right]{siunitx}
\usepackage{physics}
\usepackage{cases}

\usepackage{algorithm}
\usepackage{algpseudocode}

\usepackage[hypertexnames=false,
            colorlinks=true,
            pdfborderstyle={/S/U/W 0},
            citecolor=blue]{hyperref}
\usepackage[capitalise, noabbrev]{cleveref}

\usepackage{graphicx}
\usepackage{subcaption}

\usepackage[colorinlistoftodos,
            bordercolor=white,
            textsize=footnotesize]{todonotes}

\newfloat{algorithm}{t}{lop}

\everymath{\displaystyle}
\allowdisplaybreaks

\graphicspath{{fig}{standalone}}

\DeclareSIUnit\mmhg{mmHg}

\theoremstyle{definition}
\newtheorem{remark}{Remark}

\newcommand{\reva}[1]{#1}
\newcommand{\revb}[1]{#1}

\renewcommand{\epsilon}{\varepsilon}

\renewcommand{\tilde}{\widetilde}

\newcommand{\lifex}{\texttt{life\textsuperscript{x}}}
\newcommand{\dealii}{\texttt{deal.II}}

\newcommand{\domain}{\Omega}
\newcommand{\gammaepi}{\Gamma_\text{epi}}
\newcommand{\gammaendo}{\Gamma_\text{endo}}
\newcommand{\gammabase}{\Gamma_\text{base}}
\newcommand{\normal}{\mathbf{n}}

\newcommand{\potential}{v}
\newcommand{\ionicvars}{\mathbf{w}}
\newcommand{\actvars}{\mathbf{s}}
\newcommand{\displacement}{\mathbf{d}}
\newcommand{\circvars}{\mathbf{c}}

\newcommand{\fibers}{\mathbf{f}_0}
\newcommand{\sheets}{\mathbf{s}_0}
\newcommand{\normals}{\mathbf{n}_0}
\newcommand{\difftensor}{\mathbf{D}_\mathrm{m}}
\newcommand{\sigmaf}{\sigma_\text{f}}
\newcommand{\sigmas}{\sigma_\text{s}}
\newcommand{\sigman}{\sigma_\text{n}}
\newcommand{\Iion}{I_\text{ion}}
\newcommand{\Iapp}{I_\text{app}}
\newcommand{\calcium}{[\text{Ca}^{2+}]_\text{i}}

\newcommand{\Ftensor}{\mathbf{F}}
\newcommand{\jacobian}{J}
\newcommand{\SL}{\text{SL}}
\newcommand{\actstress}{\mathbf{P}_\text{act}}
\newcommand{\passtress}{\mathbf{P}_\text{pas}}
\newcommand{\stress}{\mathbf{P}}
\newcommand{\Kepi}{\mathbf{K}_\text{epi}}
\newcommand{\Cepi}{\mathbf{C}_\text{epi}}
\newcommand{\eye}{\mathbf{I}}

\newcommand{\plv}{p_\text{LV}}
\newcommand{\Vlv}{V_\text{LV}}

\newcommand{\src}{^\text{src}}
\newcommand{\dst}{^\text{dst}}
\newcommand{\srcpoint}[1]{\mathbf x\src_{#1}}
\newcommand{\dstpoint}[1]{\mathbf x\dst_{#1}}
\newcommand{\rbf}{\phi}
\newcommand{\rbfradius}[1]{r_{#1}}
\newcommand{\interpmatrix}{\Phi_\text{int}}
\newcommand{\evalmatrix}{\Phi_\text{eval}}
\newcommand{\cardinal}[1]{\psi_{#1}}
\newcommand{\precmatrix}{P}
\newcommand{\mesh}{\mathcal{T}}
\newcommand{\ep}{^\text{EP}}
\newcommand{\mech}{^\text{M}}
\newcommand{\mechpoint}[1]{\mathbf x\mech_{#1}}
\newcommand{\eppoint}[1]{\mathbf x\ep_{#1}}
\newcommand{\matrixU}{\mathbf{U}}
\newcommand{\matrixV}{\mathbf{V}}
\newcommand{\Sigmamatrix}{\boldsymbol\Sigma}

\newcommand{\columnU}{\mathbf{u}}
\newcommand{\columnV}{\mathbf{v}}
\newcommand{\singvalue}{\sigma}

\newcommand{\diag}{\operatorname{diag}}

\newcommand{\testvect}{\mathbf{a}}
\newcommand{\refdirection}{\mathbf{w}}

\newcommand{\interp}[1]{\Pi_{#1}}
\newcommand{\interpresc}[1]{\Pi^\text{res}_{#1}}

\newcommand{\nestedd}{nested-$\displacement$}
\newcommand{\rbfd}{RBF-$\displacement$}
\newcommand{\rbffew}{RBF-$\Ftensor$-E}
\newcommand{\rbffsvd}{RBF-$\Ftensor$-SVD}

\acsetup{make-links  = false}

\DeclareAcronym{RBF}{long={radial basis function}, short={RBF}}
\DeclareAcronym{ODE}{long={ordinary differential equation}, short={ODE}}
\DeclareAcronym{SVD}{long={singular value decomposition}, short={SVD}}
\DeclareAcronym{IMEX}{long={implicit-explicit}, short={IMEX}}
\DeclareAcronym{DoF}{long={degree of freedom}, short={DoF}, long-plural-form={degrees of freedom}}

\makeatletter
\def\ps@pprintTitle{%
 \let\@oddhead\@empty
 \let\@evenhead\@empty
 \def\@oddfoot{}%
 \let\@evenfoot\@oddfoot}
\makeatother
\biboptions{sort&compress}

\begin{document}

\begin{frontmatter}
      \title{Preserving the positivity of the deformation gradient determinant in intergrid interpolation by combining RBFs and SVD: application to cardiac electromechanics}

      \author[1]{Michele Bucelli}
      \author[1]{Francesco Regazzoni}
      \author[1]{Luca Dede'\texorpdfstring{\corref{cor1}}{}}
      \author[1,2]{Alfio Quarteroni}

      \affiliation[1]{
            organization={MOX, Laboratory of Modeling and Scientific Computing, Dipartimento di Matematica, Politecnico di Milano},
            addressline={Piazza Leonardo da Vinci 32},
            postcode={20133},
            city={Milano},
            country={Italy}}
      \affiliation[2]{
            organization={Institute of Mathematics, École Polytechnique Fédérale de Lausanne},
            addressline={Station 8, Av. Piccard},
            postcode={CH-1015},
            city={Lausanne},
            country={Switzerland (Professor Emeritus)}}

      \cortext[cor1]{Corresponding author. E-mail: luca.dede@polimi.it}
      \date{Last update: {\today}}

      \journal{}

      \begin{abstract}
            The accurate, robust and efficient transfer of the deformation gradient tensor between meshes of different resolution is crucial in cardiac electromechanics simulations. This paper presents a novel method that combines rescaled localized Radial Basis Function (RBF) interpolation with Singular Value Decomposition (SVD) to \reva{preserve the positivity} of the determinant of the deformation gradient tensor. The method involves decomposing the \reva{evaluations of the tensor at the quadrature nodes of the source mesh} into rotation matrices and diagonal matrices of singular values; computing the RBF interpolation of the quaternion representation of rotation matrices and the singular value logarithms; reassembling the deformation gradient tensors at \reva{quadrature nodes of the destination mesh, to be used in the assembly of the electrophysiology model equations}. The proposed method overcomes limitations of existing interpolation methods, including nested intergrid interpolation and RBF interpolation of the displacement field, that may lead to \reva{the loss of physical meaningfulness of the mathematical formulation and then to} solver failures \reva{at the algebraic level,} due to negative determinant values. Furthermore, the proposed method enables the transfer of solution variables between finite element spaces of different degrees and shapes and without stringent conformity requirements between different meshes, thus enhancing the flexibility and accuracy of electromechanical simulations. We show numerical results confirming that the proposed method enables \reva{the} transfer of the deformation gradient tensor, allowing to successfully run simulations in cases where existing methods fail. This work provides an efficient and \reva{robust} method for the intergrid transfer of the deformation gradient tensor, thus enabling independent tailoring of mesh discretizations to the particular characteristics of \reva{the individual physical components concurring to the} of the multiphysics model.
      \end{abstract}

\end{frontmatter}


\section{Introduction}

The multiple \reva{physical models} involved in the \reva{mathematical representation} of cardiac electromechanics \cite{niederer2019computational, piersanti20223d, quarteroni2017integrated, regazzoni2022cardiac, salvador2021electromechanical, strocchi2020simulating, stella2022fast, quarteroni2022modeling, levrero2020sensitivity, quarteroni2019mathematical, gerach2021electro, augustin2016patient, gurev2011models} and electro-fluid-mechanics \cite{bucelli2023mathematical, santiago2018fully, viola2022fsei, viola2023gpu, verzicco2022electro, zingaro2023electromechanics, karabelas2018towards} are characterized by very different spatial and temporal scales. In particular, cardiac electrophysiology \cite{collifranzone2014mathematical, piersanti2021modeling, sundnes2007computing, sung2020personalized, gillette2021framework} features very fast transients and sharp propagating fronts \cite{bucelli2021multipatch, pegolotti2019isogeometric}, requiring\reva{, upon finite element discretization,} very fine computational \reva{meshes} to be captured, with a typical size \reva{for the mesh elements} of around \SI{0.3}{\milli\metre} \cite{arevalo2016arrhythmia,niederer2011verification,trayanova2017imaging}.
\reva{Conversely}, cardiac mechanics \reva{does not require such a fine discretization} \cite{regazzoni2022cardiac,piersanti20223d}. Solving electrophysiology and mechanics using the same spatial resolution (dictated by the accuracy requirements of the former) entails an excessive computational cost. For this reason, it can be computationally convenient to solve the two problems using different spatial discretizations (a fine one for electrophysiology and a coarse one for mechanics), \reva{relying on} suitable intergrid operators to transfer variables between the two models \cite{regazzoni2022cardiac,salvador2020intergrid}. The intergrid operators must primarily transfer the intracellular calcium concentration from the electrophysiology model to the active mechanics one. Moreover, \revb{when considering} mechano-electric feedback effects \cite{collet2015numerical,salvador2022role, timmermann2017integrative, quinn2021cardiac}, the deformation gradient tensor \revb{should also be transferred} from the mechanics model \revb{back} to the electrophysiology one.

A simple approach to intergrid interpolation is based on using the same mesh for electrophysiology and mechanics, but using higher order polynomials for electrophysiology, as done e.g. in \cite{bucelli2023mathematical,fedele2023comprehensive}, where quadratic finite elements are used for electrophysiology and linear elements are used for solid mechanics. Higher polynomial orders may require the introduction of \reva{appropriate} high-order methods \cite{africa2023matrix,bucelli2021multipatch,pegolotti2019isogeometric}, which in turn may call for significant changes in existing computational pipelines.

Alternatively, the computational framework presented in \cite{piersanti20223d,regazzoni2022cardiac} relies on \reva{interpolation between nested meshes} to solve cardiac mechanics on a coarse hexahedral mesh, and electrophysiology on a finer mesh obtained by subdividing the coarse one. That approach relies on an octree implementation of the mesh data structures \cite{arndt2020dealii,burstedde2011p4est}, which makes the interpolation efficient for quadrilateral and hexahedral discretizations. Nested \reva{meshes}, however, pose significant restrictions. For example, if the solid mechanics mesh is locally refined to capture geometrical features or material discontinuities, such inhomogeneities are inherited by the electrophysiology mesh, even though a non-constant mesh resolution may lead to artificial spatial variations in conduction velocity \cite{niederer2011verification, pathmanathan2011significant}. Moreover, \reva{the geometrical detail captured by the fine mesh} is limited by that of the coarse mesh, partially countering the advantages of mesh refinement.

A more flexible approach is offered by \ac{RBF} interpolation \cite{deparis2016internodes,deparis2014rescaled,salvador2020intergrid,voet2022internodes}. In this case, the two meshes \revb{can be} independent, both geometrically and parametrically. \Ac{RBF} interpolation was applied to cardiac electromechanics in \cite{salvador2020intergrid}, showing that the interpolation allows for an accurate segregated electromechanical solver significantly faster than its monolithic counterpart \cite{gerbi2018monolithic}.

The interpolation of the deformation gradient $\Ftensor = \eye + \grad\displacement$, with $\displacement$ being the tissue displacement field, from the coarse to the fine mesh requires special care. Indeed, if mechano-electrical feedbacks are included in the model \cite{collet2015numerical, salvador2022role, timmermann2017integrative, quinn2021cardiac}, the electrophysiology equations involve the inverse of the deformation gradient $\Ftensor$ and its determinant $\jacobian$. For the problem to be well posed and physically meaningful, $\jacobian$ should be positive on the whole domain. However, naive tensor interpolation methods cannot guarantee that this is verified \cite{satheesh2022structure}.

In this work, we combine rescaled localized \ac{RBF} interpolation \cite{deparis2014rescaled} with \ac{SVD}, to obtain an interpolation method for the tensor $\Ftensor$ that preserves the sign of its determinant, in an approach similar to the one presented in \cite{satheesh2022structure}. The interpolation method is applied to ventricular electromechanical simulations, showing that it improves the robustness of the method over simpler techniques. We compare the newly introduced method against alternative interpolation methods \cite{regazzoni2022cardiac,salvador2020intergrid} in terms of numerical \revb{accuracy} and computational costs.

\revb{Our results show that the proposed method allows to successfully perform simulations where previously introduced methods would fail. Indeed, previous methods yield negative values for the deformation gradient determinant $\jacobian$ after interpolation, whereas the newly proposed one guarantees by construction its positivity. This avoids the appearence of highly unphysical values of $\jacobian$, and makes sure that the discrete problem remains well posed throughout the simulation. Furthermore, we demonstrate that the computational cost of interpolation remains negligible with respect to the overall cost of the simulation, and is linearly scalable in a parallel computing framework.}

The rest of the paper is structured as follows. \Cref{sec:electromechanics} briefly reviews the models used for cardiac electromechanics. In \cref{sec:intergrid}, we recall the rescaled localized \ac{RBF} interpolation method, and describe the interpolation procedure based on \ac{SVD}. \Cref{sec:discretization} describes the discretization strategy used for the electromechanical model, and \cref{sec:numerical-experiments} presents some numerical experiments highlighting the properties of the proposed method. Finally, in \cref{sec:conclusions}, we draw some conclusive remarks.

\section{Electromechanical modeling of the heart}
\label{sec:electromechanics}

\begin{figure}
      \centering

      \includegraphics{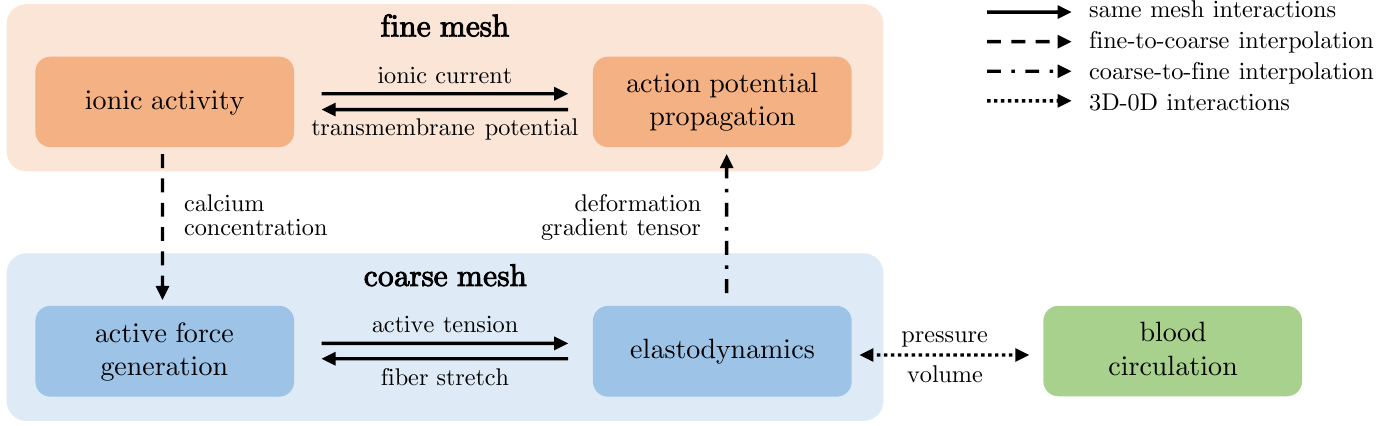}

      \caption{Cardiac electromechanical model and its subsystems. The large shaded areas enclose models that are discretized using the same computational mesh. Arrows denote the variables that realize the coupling between the subsystems. Solid lines correspond to quantities that are defined on the same computational mesh; dashed and dot-dashed lines correspond to quantities that need to be transferred from one mesh to another; the dotted line corresponds to quantities associated with the coupling between the 3D electromechanical model and the 0D model (i.e. lumped-parameter model) of blood circulation.}
      \label{fig:model}
\end{figure}

\revb{The sketch of} the \reva{electromechanics} model used in this work \revb{is displayed} in \cref{fig:model}. We refer the interested reader to \cite{regazzoni2022cardiac} for further details on \reva{the model and its} numerical solution.

\begin{figure}
      \centering

      \includegraphics{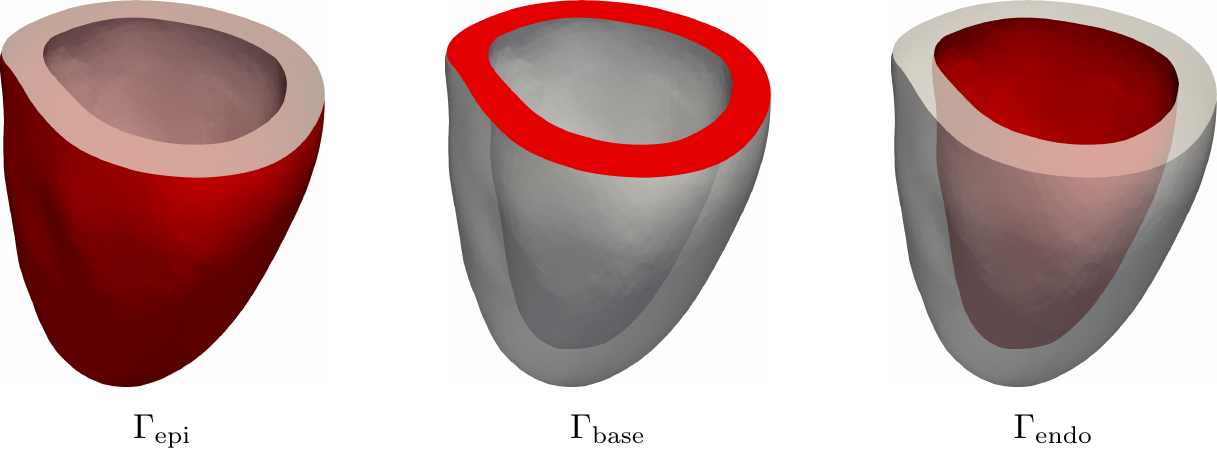}

      \caption{Computational domain $\domain$ and its boundaries $\gammaepi$, $\gammabase$ and $\gammaendo$.}
      \label{fig:domain}
\end{figure}

Let $\domain \subset \mathbb{R}^3$ be an open, bounded domain representing a left ventricle (shown in \cref{fig:domain}), and let $T > 0$. We consider a coupled problem involving cardiac electrophysiology, active force generation, muscular mechanics and the circulatory system, with the following unknowns:
\begin{align*}
      \potential &: \domain\times(0, T) \to \mathbb{R} & \text{transmembrane potential}\;, \\
      \ionicvars &: \domain\times(0, T) \to \mathbb{R}^{N_\text{ion}} & \text{ionic variables}\;, \\
      \actvars &: \domain\times(0, T) \to \mathbb{R}^{N_\text{act}} & \text{activation variables}\;, \\
      \displacement &: \domain\times(0, T) \to \mathbb{R}^3 & \text{solid displacement}\;, \\
      \circvars &: (0, T) \to \mathbb{R}^{N_\text{circ}} & \text{circulation state variables}\;.
\end{align*}

We model the evolution of the transmembrane potential $v$ with the monodomain equation \cite{collifranzone2014mathematical,sundnes2007computing}, including geometry-mediated mechano-electrical feedback effects \cite{collet2015numerical,salvador2022role,timmermann2017integrative}:
\begin{equation}\begin{cases}
      \jacobian \chi C_m \pdv{\potential}{t}
            - \div(\jacobian \Ftensor^{-1}\difftensor\Ftensor^{-T} \grad\potential)
            + \jacobian\chi\Iion(\potential, \ionicvars)
            = \jacobian\chi\Iapp & \text{in } \domain \times (0,T)\;, \\
      \jacobian\Ftensor^{-1}\difftensor\Ftensor^{-T}\grad\potential\cdot\normal = 0 & \text{on }\partial\domain \times (0, T)\;, \\
      \potential = \potential_0 & \text{in }\domain \times \{0\}\;.
\end{cases} \label{eq:monodomain} \end{equation}
In the above, $\Ftensor = \eye + \grad\displacement$ is the deformation gradient associated with the displacement of the cardiac muscle, and $\jacobian = \det\Ftensor$ is its determinant. Their presence in \eqref{eq:monodomain} accounts for the effect of the deformation onto the propagation of the electrical activation. The tensor $\difftensor$ accounts for the anisotropic conductivity of the cardiac muscle, and it is defined as
\begin{equation}
      \difftensor =
      \sigmaf\frac{\Ftensor\fibers\otimes\Ftensor\fibers}{\|\Ftensor\fibers\|^2} +
      \sigmas\frac{\Ftensor\sheets\otimes\Ftensor\sheets}{\|\Ftensor\sheets\|^2} +
      \sigman\frac{\Ftensor\normals\otimes\Ftensor\normals}{\|\Ftensor\normals\|^2}\;,
      \label{eq:difftensor}
\end{equation}
wherein $\{\fibers, \sheets, \normals\}$ is a space-dependent orthonormal triplet that represents the local direction of fibers, fiber sheets and cross-fibers, respectively \cite{piersanti2021modeling,rossi2014thermodynamically,bayer2012novel}.
The evolution of ionic variables is prescribed by a suitable ionic model, in the general form
\begin{equation}\begin{cases}
      \pdv{\ionicvars}{t} = \mathbf{F}_\text{ion}(\potential, \ionicvars) & \text{in } \domain \times (0,T)\;, \\
      \ionicvars = \ionicvars_0 & \text{in } \domain \times \{0\}\;.
\end{cases} \label{eq:ionic} \end{equation}
We consider the ventricular ionic model by ten Tusscher and Panfilov \cite{ten2006alternans}.

\begin{remark}
      Although not considered in this work, the system \eqref{eq:monodomain} can be endowed with additional reactive terms that account for other sources of mechano-electrical feedbacks, including e.g. stretch-activated currents \cite{salvador2022role}. Although they have a minor effect in sinus rhythm, they may become extremely relevant in pathological scenarios \reva{involving arrhythmias} \cite{salvador2022role}, so that an accurate evaluation of $\Ftensor$ can become extremely important.
\end{remark}

One of the ionic variables $\mathbf w$ represents the intracellular calcium concentration $\calcium$, which is used as input to the force generation model describing the evolution of the contraction state $\mathbf s$. To this end, we consider the RDQ18 \cite{regazzoni2018active}, that can be expressed as a system of \acp{ODE}:
\begin{equation}\begin{cases}
      \pdv{\actvars}{t} = \mathbf{F}_\text{act}\left(\actvars, \calcium, \SL\right) & \text{in } \domain \times (0, T)\;, \\
      \actvars = \actvars_0 & \text{in } \domain \times \{0\}\;.
\end{cases} \label{eq:force-generation} \end{equation}
In the above, SL is the sarcomere length, defined as $\SL = \SL_0 \sqrt{I_\text{4f}}$, where $I_\text{4f} = \Ftensor\fibers\cdot \Ftensor\fibers$. The contraction state is then used to compute an active stress tensor through
\begin{equation*}
      \actstress(\displacement, \actvars) = T_\text{act}^\text{max}\,P_\text{act}(\actvars)\frac{\Ftensor\fibers \otimes \Ftensor\fibers}{\sqrt{I_\text{4f}}}\;,
\end{equation*}
where $T_\text{act}^\text{max}$ is the maximum active tension and $P_\text{act}(\actvars) \in [0, 1]$ is the permissivity. We refer to \cite{regazzoni2018active} for the precise definition of $\mathbf{F}_\text{act}$ and $P_\text{act}$. Due to the complexity of the model, we consider a surrogate version obtained through artificial neural networks, as presented in \cite{regazzoni2019machine,regazzoni2020machine}.

The displacement of the muscle is modeled by the elastodynamics equation in the hyperelastic framework:
\begin{equation}\begin{cases}
      \rho\pdv[2]{\displacement}{t} - \div \stress(\displacement, \actvars) = \mathbf 0 & \text{in } \domain \times (0, T)\;, \\
      \stress(\displacement, \actvars)\normal + \Kepi \displacement + \Cepi \pdv{\displacement}{t} = \mathbf 0 & \text{on } \gammaepi \times (0, T)\;, \\
      \stress(\displacement, \actvars)\normal = \plv(t) |\jacobian\Ftensor^{-T}\normal| \frac{\int_{\gammaendo}\jacobian\Ftensor^{-T}\normal\,d\Gamma}{\int_{\gammabase}|\jacobian\Ftensor^{-T}\normal|\,d\Gamma} & \text{on } \gammabase \times (0, T)\;, \\
      \stress(\displacement, \actvars)\normal = -\plv(t) \jacobian\Ftensor^{-T}\normal & \text{on } \gammaendo\times (0, T)\;,
\end{cases} \label{eq:mechanics} \end{equation}
where $\stress(\displacement, \actvars) = \actstress(\displacement, \actvars) + \passtress(\displacement)$ is the first Piola-Kirchhoff stress tensor, decomposed into the sum of the active stress $\actstress$ and the passive stress $\passtress$. The latter is obtained as the derivative of a suitable strain energy density function representing the constitutive law of the solid:
\begin{equation*}
      \passtress = \pdv{\mathcal W}{\Ftensor}\;.
\end{equation*}
We consider the orthotropic constitutive law of \cite{guccione1991finite,usyk2002computational}, with an additional penalization term imposing near-incompressibility as presented in \cite{regazzoni2022cardiac}. $\gammaepi$, $\gammabase$ and $\gammaendo$ are subsets of $\partial\domain$, shown in \cref{fig:domain}, and $\normal$ is the outward-directed normal unit vector. The boundary condition on $\gammaepi$ expresses the interaction of the heart with the pericardium and the surrounding tissue \cite{pfaller2019importance, strocchi2020simulating}. The tensors $\Kepi$ and $\Cepi$ are defined as
\begin{align*}
      \Kepi = K_\text{epi}^\perp (\normal\otimes\normal)
      + K_\text{epi}^\parallel (\eye - \normal\otimes\normal)\;, \\
      \Cepi = C_\text{epi}^\perp (\normal\otimes\normal)
      + C_\text{epi}^\parallel (\eye - \normal\otimes\normal)\;.
\end{align*}
The condition on $\gammabase$ is known as energy-consistent boundary condition, and it surrogates the portion of ventricle that has been cut from the computational domain at the base \cite{regazzoni2020machine, regazzoni2022cardiac}.

The intracardiac pressure $\plv$ is obtained by solving a lumped-parameter model for the circulatory system \cite{blanco20103d,hirschvogel2017monolithic,regazzoni2022cardiac}. The model can be expressed as a system of algebraic-differential equations:
\begin{equation}\begin{cases}
      \mathbf{F}_\text{circ}\left(\dv{\circvars}{t}, \circvars, t\right) = \mathbf 0 & \text{in }(0, T)\;, \\
      \circvars(0) = \mathbf c_0\;.
\end{cases} \label{eq:circulation} \end{equation}
Two of the variables of the circulation state $\circvars$ are the left ventricular pressure $\plv$ and the left ventricular volume $\Vlv$. The coupling of the mechanics equations \eqref{eq:mechanics} and of the circulation system \eqref{eq:circulation} is obtained by imposing the pressure condition on $\gammaendo$ in \eqref{eq:mechanics}, together with the constraint $\Vlv = \Vlv^\text{3D}$, where $\Vlv^\text{3D}$ is the volume enclosed by $\gammaendo$ in the deformed configuration, computed as described in \cite{regazzoni2022cardiac}.

\section{Intergrid interpolation for electromechanics}
\label{sec:intergrid}

In the following sections, we recall the \ac{RBF} interpolation method, and present the algorithm used in this work for interpolating the deformation gradient $\Ftensor$ from the mesh of the mechanics problem to that of the electrophysiology problem.

\begin{figure}
      \centering

      \begin{subfigure}{6cm}
            \centering
            \includegraphics{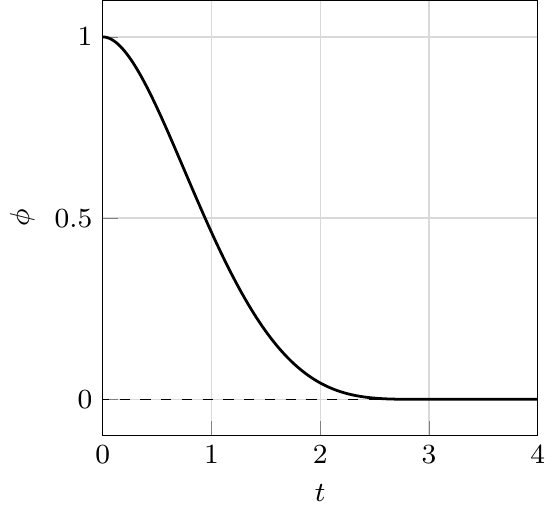}
            \caption{}
            \label{fig:rbf}
      \end{subfigure}
      \begin{subfigure}{6cm}
            \centering
            \includegraphics{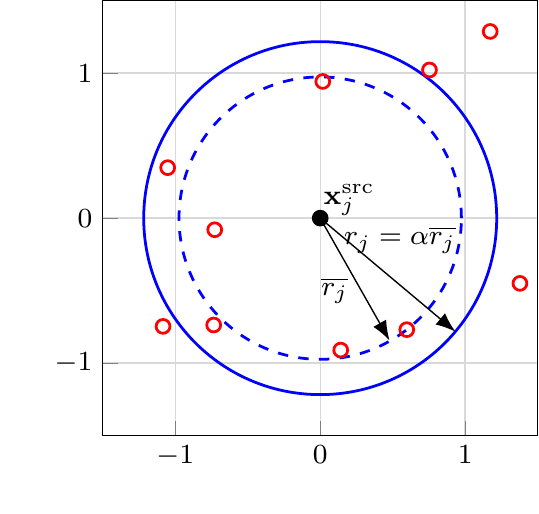}
            \caption{}
            \label{fig:adaptive-radius}
      \end{subfigure}

      \caption{(a) Plot of the $C^2$ Wendland basis function $\rbf(t, 3)$. (b) Selection of the adaptive \ac{RBF} radius $r_j$ for point $\srcpoint{j}$. The radius is such that the support encloses $M = \num{4}$ interpolation points (indicated by the red void circles), \revb{with} $\alpha = \num{1.2}$ \revb{(see \eqref{eq:rbf-radius})}.}
\end{figure}

\subsection{Rescaled localized radial basis function interpolation}
\label{sec:rescaled-localized-rbf}

Let $\{\srcpoint{i}\}_{i = 1}^{N\src}$, with $\srcpoint{i} \in \mathbb{R}^3$, be a set of distinct points. Let $f:\mathbb{R}^3 \to \mathbb{R}$ be a function, and let $f\src_i = f(\srcpoint{i})$. The \ac{RBF} interpolant of $f$ at the points $\srcpoint{i}$ is a function $\interp{f}:\mathbb{R}^3 \to \mathbb{R}$ in the form
\begin{equation}
      \interp{f}(\mathbf x) = \sum_{j = 1}^{N\src} \gamma_j\,\rbf\!\left(\|\mathbf x - \srcpoint{j}\|, \rbfradius{j}\right)\;,
      \label{eq:rbf-interpolant}
\end{equation}
where $\gamma_j$ are the interpolation coefficients, $\rbf$ is the \ac{RBF} and $\rbfradius{j}$ is the \ac{RBF} support radius associated with each point $\srcpoint{j}$. Following \cite{deparis2014rescaled,salvador2020intergrid}, we consider the compactly supported $C^2$ Wendland basis function \cite{wendland1995piecewise}, shown in \cref{fig:rbf} and defined as
\begin{equation*}
      \rbf(t, r) = \max\left\{1 - \frac{t}{r}, 0\right\}^4 \left(1 + 4\frac{t}{r}\right)\;.
\end{equation*}
The coefficients $\gamma_i$ in \eqref{eq:rbf-interpolant} are determined by imposing the interpolation condition $\interp{f}(\srcpoint{i}) = f\src_i$ for all $i = 1, 2, \dots, N\src$. This gives rise to the linear system
\begin{equation}
      \interpmatrix \boldsymbol\gamma = \boldsymbol f\src\;,
      \label{eq:rbf-linear-system}
\end{equation}
where $\interpmatrix \in \mathbb{R}^{N\src \times N\src}$ is a matrix whose entries are
\begin{equation*}
      (\interpmatrix)_{ij} = \rbf(\|\srcpoint{i} - \srcpoint{j}\|, \rbfradius{j})\;,
\end{equation*}
with $\boldsymbol\gamma = (\gamma_1, \gamma_2, \dots, \gamma_{N\src})^T \in \mathbb{R}^{N\src}$ and $\mathbf f\src = (f\src_1, f\src_2, \dots, f\src_{N\src})^T \in \mathbb{R}^{N\src}$.

Assume now that the interpolant must be evaluated on a new set of points $\{\dstpoint{i}\}_{i = 1}^{N\dst}$. Let $f\dst_i = \interp{f}(\dstpoint{i})$. There holds, for all $i = 1, 2, \dots, N\dst$:
\begin{equation*}
      f\dst_i = \sum_{j = 1}^{N\src} \gamma_j \rbf(\|\dstpoint{i} - \srcpoint{j}\|, \rbfradius{j})\;,
\end{equation*}
which can be expressed compactly as
\begin{equation*}
      \mathbf f\dst = \evalmatrix \boldsymbol\gamma\;,
\end{equation*}
where $\evalmatrix \in \mathbb{R}^{N\dst \times N\src}$ is a matrix whose entries are
\begin{equation*}
      (\evalmatrix)_{ij} = \phi(\|\dstpoint{i} - \srcpoint{j}\|, \rbfradius{j})\;,
\end{equation*}
and $\mathbf f\dst = (f\dst_1, f\dst_2, \dots, f\dst_{N\dst})^T \in \mathbb{R}^{N\dst}$.

The interpolation radii $\rbfradius{j}$ are selected adaptively. Intuitively, the radius should be smaller in regions of space where the interpolation points $\srcpoint{j}$ are more densely clustered. In \cite{deparis2014rescaled}, the authors select the radius based on the connectivity of the mesh where the interpolated data is defined. In a parallel computing context, where the mesh is distributed across several processors, the connectivity may not be easily accessible. Moreover, this notion does not generalize to the case where the interpolation and evaluation points are not the nodes of a mesh (as considered in \cref{sec:interpolation-quadrature}). For this reason, we follow an approach similar to the one used in \cite{voet2022internodes}, and choose
\begin{equation}
      \rbfradius{j} = \alpha \overline{\rbfradius{j}}\;,
      \label{eq:rbf-radius}
\end{equation}
where $\overline{\rbfradius{j}}$ is the smallest value such that the sphere centered at $\srcpoint{j}$ with radius $\overline{\rbfradius{j}}$ contains at least $M$ other interpolation points. This procedure is represented in \cref{fig:adaptive-radius}.

As discussed in \cite{deparis2014rescaled}, the \ac{RBF} interpolation procedure described above may yield large oscillations in the interpolant, and is very sensitive to the choice of the interpolation radius. To avoid this issue, \reva{we consider} the following rescaled interpolant\reva{, introduced in \cite{deparis2014rescaled}}:
\begin{equation*}
      \interpresc{f}(\mathbf x) = \frac{\interp{f}(\mathbf x)}{\interp{g}(\mathbf x)}\;,
\end{equation*}
where $\interp{g}$ is the \ac{RBF} interpolant of the constant function $g(\mathbf x) = 1$ at the nodes $\srcpoint{i}$. All the results presented in this paper make use of the rescaled interpolant.


The interpolation of vector fields is obtained by separately interpolating each component.

\subsection{Preconditioning of the interpolation system by means of approximate cardinal functions}
\label{sec:preconditioner}

The matrix $\interpmatrix$ is sparse due to the \ac{RBF} $\rbf$ having compact support. We solve system \eqref{eq:rbf-linear-system} by means of the preconditioned GMRES method \cite{saad2003iterative}. The choice of a suitable preconditioner is crucial towards the efficient construction of the interpolant.

Following \cite{beatson1999fast,brown2005approximate,gumerov2007fast}, we consider a preconditioner based on approximated cardinal functions. The procedure to assemble it is described in \cref{alg:preconditioner}. The idea behind the preconditioner is to construct, for each point $\srcpoint{i}$, an approximate cardinal function $\cardinal{i}$, such that
\begin{align*}
      \cardinal{i}(\mathbf x) &= \sum_{j \in S^i} \lambda^i_j\,\rbf(\|\srcpoint{j} - \mathbf x\|, \rbfradius{j})\;, \\
      \cardinal{i}(\srcpoint{i}) &= 1\;, \\
      \cardinal{i}(\srcpoint{j}) &= 0 \quad \text{for } j \in S^i,\; j \neq i\;,
\end{align*}
where $S^i$ is a set of points that are sufficiently close to $\srcpoint{i}$. Finding the cardinal function is itself an interpolation problem over the points specified by $S^i$. The preconditioner matrix $\precmatrix^{-1} \in \mathbb{R}^{N\src \times N\src}$ has entries
\begin{equation*}
      \left(\precmatrix^{-1}\right)_{ji} = \begin{cases}
            \lambda^i_j & \text{if } j \in S^i\,, \\
            0 & \text{otherwise.}
      \end{cases}
\end{equation*}
We select $S^i$ to be the set of indices $j$ such that $(\interpmatrix)_{ij} \neq 0$. We refer to \cite{beatson1999fast,brown2005approximate,gumerov2007fast} for further considerations on this preconditioning strategy.

\begin{algorithm}
      \caption{Construction of the preconditioner $\precmatrix$ for system \eqref{eq:rbf-linear-system}.}
      \label{alg:preconditioner}

      \begin{algorithmic}[1]
        \Require interpolation matrix $\interpmatrix$
        \Ensure preconditioner matrix $\precmatrix$
        \vspace{\baselineskip}

        \State initialize $\precmatrix = 0$
        \For{$i = 0, 1, \dots, N\src$}
          \State $S^i \gets \{ j \in \mathbb{N} : (\interpmatrix)_{ij} \neq 0 \}$, and let $s^i_m$ be its elements, with $m = 0, 1, \dots, n_i$

          \State 

          \State compute the matrix $L^i \in \mathbb{R}^{n_i \times n_i}$\reva{, with} entries $L^i_{lm} = \rbf(\|\srcpoint{s_l} - \srcpoint{s_m}\|, r_{s_m})$
          \State compute the vector $\mathbf{r}^i \in \mathbb{R}^{n_i}$\reva{, with} entries $\mathbf{r}^i_l = \delta_{i s^i_l}$
          \State solve the system $L^i \boldsymbol{\lambda}^i = \mathbf{r}^i$ with the GMRES method

          \State 

          \For{$m = 0, 1, \dots, n_i$}
            \State $\left(\precmatrix^{-1}\right)_{s^i_m,i} \gets \boldsymbol{\lambda}^i_m$
          \EndFor
        \EndFor
      \end{algorithmic}
\end{algorithm}

The construction of the preconditioner \reva{can be} very costly, requiring the solution of a small yet dense linear system for each point $\srcpoint{i}$. However, if the interpolant must be constructed and evaluated several times, the computational saving in the repeated solution of \eqref{eq:rbf-linear-system} amortizes the cost of the preconditioner assembly.

We use a Gauss-Seidel preconditioner to accelerate the convergence of the GMRES method for the linear system $L^i \boldsymbol{\lambda}^i = \mathbf{r}^i$ in \cref{alg:preconditioner}. Due to the approximate nature of the preconditioner, the interpolation problem to compute $\cardinal{i}$ can be solved to low accuracy without significantly damaging the effectiveness of $\precmatrix$. In all the tests described in this paper, we use a relative tolerance of \num{1e-1}.

\begin{remark}
      With respect to the strategy outlined in \cite{beatson1999fast}, our choice of $S^i$ allows to compute $\precmatrix$ by efficiently reusing data structures and quantities already computed for the construction of $\interpmatrix$. We remark that this choice is not the only possible one, and alternative options may result in a better preconditioner, at the price of additional computational overhead for its initialization.
\end{remark}

\subsection{RBF interpolation in the finite element framework}

In the context of finite elements, let us consider a domain $\Omega \subset \mathbb{R}^3$, and two independent meshes approximating $\Omega$, $\mesh\src$ and $\mesh\dst$. Each of the two meshes can be composed of tetrahedral or hexahedral elements (not necessarily the same on both meshes). On the two meshes, we consider the finite element spaces $V\src$ and $V\dst$, composed of piecewise polynomials of degree $p\src$ and $p\dst$, respectively, spanned by suitable interpolatory basis functions (such as the Lagrangian basis functions \cite{hughes2012finite,quarteroni2017numerical}) on meshes $\mesh\src$ and $\mesh\dst$.

We consider two different interpolation strategies: interpolation between \acp{DoF} and interpolation between quadrature nodes. Alternative strategies (e.g. interpolation from \acp{DoF} to quadrature nodes, or vice versa) can also be considered, although they are not relevant for our target application.

\subsubsection{Interpolation between degrees of freedom}
\label{sec:interpolation-dofs}

The points $\srcpoint{i}$ are the support points of the \acp{DoF} on $\mesh\src$, and the points $\dstpoint{i}$ are those on $\mesh\dst$. If $f$ is a function belonging to the finite element space $V\src$, then the vector $\mathbf f\src$ is the vector of its control variables. Similarly, $\mathbf f\dst$ is the vector of control variables of the finite element interpolation of $\interp{f}^\text{res}$ in $V\dst$. Therefore, the computation of the interpolant as described in \cref{sec:rescaled-localized-rbf} allows to easily interpolate from the finite element space $V\src$ onto the space $V\dst$. This is the approach we follow to interpolate the intracellular calcium concentration $\calcium$ from the mesh used for electrophysiology onto the one used for mechanics.

\subsubsection{Interpolation between quadrature nodes}
\label{sec:interpolation-quadrature}

In some instances, one may need to interpolate a \reva{function} that is not well defined at the \acp{DoF} of $\mesh\src$. That is the case of the interpolation of $\Ftensor$ from the mechanics to the electrophysiology mesh. Indeed, since the finite element space is globally continuous but only piecewise differentiable, $\Ftensor$ is not well defined at those \acp{DoF} that lie on the boundary of a mesh element.

In that case, we choose the points $\srcpoint{i}$ to be internal to the mesh elements, by considering suitable Gaussian quadrature nodes for each element. We consider a Gaussian quadrature formula with $q$ quadrature nodes in each coordinate direction.

As typical in finite elements, we approximate the integrals arising from the Galerkin formulation of \eqref{eq:monodomain} by using Gaussian quadrature. Therefore, we need to evaluate $\Ftensor$ on the quadrature nodes of the elements of $\mesh\dst$, and it is convenient to select $\dstpoint{i}$ to be those quadrature nodes.

\begin{remark}
      The source points $\srcpoint{i}$ do not need to be nodes of a quadrature formula, and any points lying on the interior of elements can be used. However, Gaussian quadrature points are usually convenient from an implementation viewpoint, since they can be easily computed in any finite element software, and the evaluaton of $\Ftensor$ at those point is easily accessible.
\end{remark}

\subsection{Interpolation of the deformation gradient}
\label{sec:interpolation-deformation-gradient}

For \eqref{eq:monodomain} to be well-posed after discretization, the interpolation of $\Ftensor$ should be done in a way that preserves the \reva{positivity} of $\jacobian = \det\Ftensor > 0$: indeed, if negative values of $\jacobian$ arise in \eqref{eq:monodomain}, its numerical solution may diverge.
However, the set of tensors with positive determinant is not a linear space (hence, linear combinations of positive-determinant tensors might have negative determinant), nor it is convex (not even convex combinations of positive-determinant tensors are guaranteed to preserve the sign of the determinant).
This explains why a naive interpolation of the deformation gradient tensor might yield non-physical results.

In \cite{salvador2020intergrid}, the deformation gradient is evaluated \reva{at the quadrature points of} the fine mesh by combining the interpolation between \acp{DoF} of the displacement with the Zienkiewicz-Zhu gradient recovery technique \cite{zienkiewicz1992superconvergent}. This approach, however, does not guarantee $\jacobian > 0$ after interpolation, and indeed we observed that in some situations it results in a breakdown of the numerical solver. To enforce $\jacobian > 0$, we combine \ac{RBF} interpolation between quadrature nodes (\cref{sec:interpolation-quadrature}) with \ac{SVD}, in an approach similar to \cite{satheesh2022structure}.

Let $\mesh\mech$ be the mesh used for the mechanics problem, and $\mesh\ep$ the one used for the electrophysiology problem. Let $\displacement$ be the displacement field, defined on $\mesh\mech$. We denote by $\mechpoint{i}$, with $i = 0, 1, \dots, N\mech_\text{q}$, the quadrature nodes on mesh $\mesh\mech$, and by $\eppoint{j}$, with $j = 0, 1, \dots, N\ep_\text{q}$, those on mesh $\mesh\ep$. The procedure to evaluate the interpolation of $\Ftensor$ onto the quadrature nodes $\eppoint{j}$, denoted by $\Ftensor\ep$, consists of decomposing the tensors to be interpolated into simpler objects, interpolating the latter, and finally recomposing the tensor in the quadrature nodes $\eppoint{j}$. The advantage is that such objects belong to spaces with a more suitable structure for interpolation than the space of positive-determinant tensors.

The procedure, reported in \cref{alg:rbf-svd}, assumes that the following routines exist:
\begin{itemize}
      \item \texttt{rotationToQuaternion}, to convert rotation matrices into their quaternion representation \cite{shoemake1985animating};
      \item \texttt{quaternionToRotation}, to convert a quaternion into its rotation matrix representation (normalizing it if necessary) \cite{shoemake1985animating};
      \item \texttt{SVD}, to compute the \ac{SVD} of a given input matrix, with the alignment procedure of \cref{sec:principal-direction-aligment procedure}.
\end{itemize}
\Cref{alg:rbf-svd} starts by computing the \ac{SVD} factorization of the deformation gradient $\Ftensor_i = \Ftensor(\mechpoint{i})$ at each point $\mechpoint{i}$, thus expressing it as
\begin{equation*}
      \Ftensor_i = \matrixU_i \Sigmamatrix_i \matrixV_i^T\;,
\end{equation*}
where $\matrixU_i$ and $\matrixV_i$ are rotation matrices, and $\Sigmamatrix_i = \diag(\singvalue_i^1, \singvalue_i^2, \singvalue_i^3)$ is a diagonal matrix whose diagonal entries are the singular values of $\Ftensor_i$. The matrices $\matrixU_i$ and $\matrixV_i$ are converted to the corresponding quaternions $\mathbf{q}_{\matrixU_i}$ and $\mathbf{q}_{\matrixV_i}$, respectively, and \ac{RBF} interpolation is used to evaluate at the points $\eppoint{j}$ the resulting \num{11} scalar fields (the three singular values and the four quaternion components of $\mathbf{q}_{\matrixU_i}$ and $\mathbf{q}_{\matrixV_i}$). Then, the deformation gradient $\Ftensor\ep$ is reconstructed at every point $\eppoint{j}$, by converting the quaternions back to rotation matrices and reassembling the \ac{SVD} factors. To guarantee that the determinant of $\Ftensor\ep$ remains positive, we interpolate the logarithm of the singular values, and then take their exponential after interpolation \cite{satheesh2022structure}. \revb{Indeed, there holds:}
\begin{equation*}
      \revb{
      \det\Ftensor\ep_j = \det\Sigmamatrix\ep_j = 
      \exp(\interpresc{\log\singvalue^{1}}(\eppoint{j}))\,\exp(\interpresc{\log\singvalue^{2}}(\eppoint{j}))\,\exp(\interpresc{\log\singvalue^{3}}(\eppoint{j})) > 0\;.
      }
\end{equation*}
The \ac{SVD} of the tensor $\Ftensor_i$ is not unique. The procedure we follow to define a unique decomposition is described in \cref{sec:principal-direction-aligment procedure}.

\begin{algorithm}
      \caption{Interpolation of the deformation gradient $\Ftensor$ from the mesh $\mesh\mech$ to the mesh $\mesh\ep$.}
      \label{alg:rbf-svd}

      \begin{algorithmic}[1]
            \Require evaluations of $\grad\displacement$ \reva{at} quadrature nodes $\mechpoint{i}$
            \Ensure evaluations of $\Ftensor\ep$ \reva{at} quadrature nodes $\eppoint{j}$
            \vspace{\baselineskip}

            \For{$i = 0, 1, \dots, N\mech_\text{q}$}
            \State $\Ftensor_i = \Ftensor(\mechpoint{i}) \gets I + \grad\displacement(\mechpoint{i})$
            \State $(\matrixU_i, \Sigmamatrix_i, \matrixV_i^T) \gets \texttt{SVD}(\Ftensor_i)$
            \Comment{With the alignment procedure of \cref{sec:principal-direction-aligment procedure}.}
            \State $(a_{\matrixU_i}, b_{\matrixU_i}, c_{\matrixU_i}, d_{\matrixU_i}) \gets \texttt{rotationToQuaternion}(\matrixU_i)$
            \State $(a_{\matrixV_i}, b_{\matrixV_i}, c_{\matrixV_i}, d_{\matrixV_i}) \gets \texttt{rotationToQuaternion}(\matrixV_i)$
            \EndFor

            \State 

            \State build the interpolants $\interpresc{a_{\matrixU}}$, $\interpresc{b_{\matrixU}}$, $\interpresc{c_{\matrixU}}$, $\interpresc{d_{\matrixU}}$, $\interpresc{a_{V}}$, $\interpresc{b_{V}}$, $\interpresc{c_{V}}$, $\interpresc{d_{V}}$, $\interpresc{\log\singvalue^{1}}$, $\interpresc{\log\singvalue^{2}}$, $\interpresc{\log\singvalue^{3}}$

            \State 

            \For{$j = 0, 1, \dots, N\mech_\text{q}$}
            \State $\matrixU\ep_j \gets \texttt{quaternionToRotation}\left(\interpresc{a_{\matrixU}}(\eppoint{j}), \interpresc{b_{\matrixU}}(\eppoint{j}), \interpresc{c_{\matrixU}}(\eppoint{j}), \interpresc{d_{\matrixU}}(\eppoint{j})\right)$
            \State $\matrixV\ep_j \gets \texttt{quaternionToRotation}\left(\interpresc{a_{V}}(\eppoint{j}), \interpresc{b_{V}}(\eppoint{j}), \interpresc{c_{V}}(\eppoint{j}), \interpresc{d_{V}}(\eppoint{j})\right)$
            \State $\Sigmamatrix\ep_j \gets \diag\left(\exp(\interpresc{\log\singvalue^{1}}(\eppoint{j})), \exp(\interpresc{\log\singvalue^{2}}(\eppoint{j})), \exp(\interpresc{\log\singvalue^{3}}(\eppoint{j}))\right)$
            \State $\Ftensor\ep_j = \Ftensor\ep(\eppoint{j}) \gets (\matrixU\ep)_j^T \Sigmamatrix\ep_j \matrixV\ep_j$
            \EndFor
      \end{algorithmic}
\end{algorithm}

\begin{remark}
      Both the matrices $\interpmatrix$ and $\evalmatrix$, involved in the construction and evaluation of the interpolants of \cref{alg:rbf-svd}, only depend on the location of the points $\mechpoint{i}$ and $\eppoint{j}$, which are the same for all interpolants and do not change over time. Therefore, the matrices are computed only once during the initialization phase.
\end{remark}

\begin{remark}
      The procedure described in \cref{alg:rbf-svd} performs linear interpolation between quaternions, as opposed to spherical interpolation \cite{shoemake1985animating}, allowing the whole interpolation procedure to be linear.
\end{remark}

\subsection{Procedure for the alignment of singular vectors}
\label{sec:principal-direction-aligment procedure}

In this section we describe the procedure (reported in \cref{alg:svd-alignment}) we propose in order to deal with the indeterminacy of the \ac{SVD} of the tensors to be interpolated.
We notice that the \ac{SVD} of a generic \reva{second-order tensor} $\Ftensor_i$ can be written as follows 
\begin{equation} 
      \Ftensor_i 
      = \matrixU_i \Sigmamatrix_i \matrixV_i^T
      = \sum_{j=1}^3 \singvalue_i^j \columnU_i^j \otimes \columnV_i^j \;,
      \label{eq:SVD_outer_products}
\end{equation}
where $\columnU_i^j$ (respectively, $\columnV_i^j$) is the $j$-th column of $\matrixU_i$ (respectively, $\matrixV_i$), known as left (respectively, right) singular vector.
From \eqref{eq:SVD_outer_products} it is clear that the decomposition of $\Ftensor_i$ is unaffected by (i) simultaneous reordering of the singular values and of the columns of $\matrixU_i$ and $\matrixV_i$, and by (ii) change of sign of corresponding columns of $\matrixU_i$ and $\matrixV_i$.
Actually, (i) and (ii) are the only sources of indeterminacy \cite{roman2005advanced}.
Conventionally, singular values are ordered in a decreasing manner, consequently defining the ordering of the columns of $\matrixU_i$ and $\matrixV_i$.
In this work, however, we adopt a different strategy for (i) ordering singular values/vectors and for (ii) defining the orientation of singular vectors.

\begin{algorithm}
      \caption{Reordering procedure for singular values and vectors}
      \label{alg:svd-alignment}

      \begin{algorithmic}[1]
        \Require singular values and vectors $\tilde\singvalue_i$, $\tilde\columnV_i$, $\tilde\columnU_i$, for $i = 1, 2, 3$, in order of decreasing singular values
        \Require reference triplet $(\refdirection^1, \refdirection^2, \refdirection^3)$
        \Ensure reordered singular values and vectors $\singvalue_i$, $\columnV_i$, $\columnU_i$, for $i = 1, 2, 3$
        \vspace{\baselineskip}

        \State $j_1 \gets \text{arg}\min_{i}|\refdirection^1 \cdot \tilde{\columnV}_i|$
        \State $\singvalue_1 \gets \tilde\singvalue_{j_1}$
        \If{$\left(\refdirection^1 \cdot \columnV_{j_1}\right) < 0$}
          \State $\columnV_1 \gets -\tilde\columnV_{j_1}, \quad \columnU_1 \gets -\tilde\columnU_{j_1}$
          \Else
          \State $\columnV_1 \gets \tilde\columnV_{j_1}, \quad \columnU_1 \gets \tilde\columnU_{j_1}$
        \EndIf

        \State 

        \State $j_2 \gets \text{arg}\min_{i \neq j_1}|\refdirection^2 \cdot \tilde\columnV_i|$
        \State $\singvalue_2 \gets \tilde\singvalue_{j_2}$
        \If{$\left(\refdirection^2 \cdot \tilde\columnV_{j_2}\right) < 0$}
          \State $\columnV_2 \gets -\tilde\columnV_{j_2}, \quad \columnU_2 \gets -\tilde\columnU_{j_2}$
        \Else
          \State $\columnV_2 \gets \tilde\columnV_{j_2}, \quad \columnU_2 \gets \tilde\columnU_{j_2}$
        \EndIf

        \State 

        \State $j_3 \gets \{1, 2, 3\} \backslash \{j_1, j_2\}$
        \State $\singvalue_3 \gets \tilde\singvalue_{j_3}$
        \If{$\det(\columnV_1, \columnV_2, \tilde\columnV_{j_3}) < 0$}
          \State $\columnV_3 \gets -\tilde\columnV_{j_3}, \quad \columnU_3 \gets -\tilde\columnU_{j_3}$
        \Else
          \State $\columnV_3 \gets \tilde\columnV_{j_3}, \quad \columnU_3 \gets \tilde\columnU_{j_3}$
        \EndIf
      \end{algorithmic}
\end{algorithm}

To \revb{better explain our} procedure, let us consider the application of a tensor $\Ftensor_i$ to a generic test vector $\testvect$:
\begin{equation} 
      \Ftensor_i \testvect
      = \matrixU_i \Sigmamatrix_i \matrixV_i^T \testvect
      = \sum_{j=1}^3 \singvalue_i^j (\columnV_i^j \cdot \testvect) \columnU_i^j \;.
      \label{eq:SVD_application}
\end{equation}
The application of $\Ftensor_i$ to the vector $\testvect$ can thus be interpreted as a three-step procedure: (i) we compute the components of $\testvect$ with respect to the orthonormal basis $(\columnV_i^1, \columnV_i^2, \columnV_i^3)$; (ii) we rescale them by the corresponding singular values; (iii) we compute a linear combination of the left singular vectors with coefficients computed in the previous step.
The ordering of the terms at right-hand side of \eqref{eq:SVD_application} clearly does not impact the result of the sum.
However, it affects the results of the interpolation, as the singular values $\singvalue_i^j$ (more precisely, their logarithm) with the same index $j$ are interpolated among points.
Since, in \eqref{eq:SVD_application}, the $j$-th singular value plays the role of rescaling the product $(\columnV_i^j \cdot \testvect)$, in this work we create a correspondence among different points $\mechpoint{i}$ (with $i = 0, 1, \dots, N\mech_\text{q}$) according to the directions of the right singular vectors $\columnV_i^j$.
More precisely, we match singular values that correspond to right singular vectors that are most closely aligned with each other.

With this aim, we define a reference orthonormal ordered triplet (which can be the canonical basis of Euclidean space, or the triplet that locally defines the direction of the fibers and sheets), and perform a reordering of the singular values and vectors, as well as a re-orientation of the latter, so as to maximize the alignment of the right singular vectors with the reference triplet.
In this way, when we interpolate the singular values associated with different points $\mechpoint{i}$ (with $i = 0, 1, \dots, N\mech_\text{q}$), we match singular values corresponding to right singular vectors that are as aligned as possible.


More in detail, the alignment procedure is the following.
Let us denote by $(\refdirection_i^1, \refdirection_i^2, \refdirection_i^3)$ the reference triplet, possibly depending on the point $\mechpoint{i}$.
We select as first column of the matrix $\matrixV_i$ the right singular vector that maximizes the quantity $|\refdirection_i^1 \cdot \columnV_i^j|$ over $j = 1,2,3$.
If \reva{$\refdirection_i^1 \cdot \columnV_i^j < 0$, that is if the cosine of the angle between $\refdirection_i^1$ and $\columnV_i^j$ is negative}, then we invert the sign of the components of the singular vector.
Then, we select as second column of the matrix $\matrixV_i$ the maximizer, between the remaining two right singular vectors, of $|\refdirection_i^2 \cdot \columnV_i^j|$.
As before, if $\refdirection_i^2 \cdot \columnV_i^j < 0$ we invert the sign of its components.
The last column of $\matrixV_i$ is assigned as the remaining singular vector, and its sign is selected so that $\det(\matrixV_i) > 0$.
Finally, we reorder the singular values and the left singular vectors, and we change the sign of the latter consistently with what done for the right singular vectors.

In this work, we consider the reference triplet to be the canonical basis of the Euclidean space, i.e. $\refdirection_i^k = \mathbf{e}^k$ for $k \in \{1, 2, 3\}$ and all $i \in \{0, 1, \dots, N\mech_\text{q}\}$.

\section{Numerical \texorpdfstring{\reva{approximation}}{approximation}}
\label{sec:discretization}

We employ the finite element method \cite{hughes2012finite,quarteroni2017numerical} for the spatial discretization of \eqref{eq:monodomain} and \eqref{eq:mechanics}. We consider either tetrahedral or hexahedral grids, with finite elements of order either \num{1} or \num{2} (although the methods proposed here generalize naturally to elements of higher order).

We discretize in time using finite difference schemes of order \num{1}. The monodomain equation \eqref{eq:monodomain} is discretized in a semi-implicit way, by treating explicitly the ionic current term, thus resulting in a linear problem. The ionic model is discretized with an \ac{IMEX} scheme, allowing for its direct solution \cite{regazzoni2022cardiac}. The mechanics model \eqref{eq:mechanics} is discretized with an implicit formulation, and the circulation model \eqref{eq:circulation} is solved with an explicit Euler scheme.

The \reva{electrophysiology, force generation, mechanics and circulation models} are coupled in a segregated staggered way\reva{, that is they are solved independently and sequentially at every time step \cite{regazzoni2022cardiac}}. The electrophysiology model is solved with a finer temporal discretization than the other models, to satisfy its stricter accuracy requirements \cite{fedele2023comprehensive,regazzoni2022cardiac}. Let $\Delta t$ be the time discretization step, and $\Delta t\ep = \Delta t / n\ep$ be a smaller discretization step used for mechanics. Then, the time advancing scheme is the following:
\begin{enumerate}
  \item solve $n\ep$ time steps of the electrophysiology model \eqref{eq:monodomain} and \eqref{eq:ionic}, with time step $\Delta t\ep$; \label{it:step-ep}
  \item interpolate the intracellular calcium concentration from $\mesh\ep$ to $\mesh\mech$ using interpolation between \acp{DoF} (\cref{sec:interpolation-dofs}); \label{it:interpolation-calcium}
  \item solve the force generation model \eqref{eq:force-generation} using time discretization step $\Delta t$; \label{it:step-force-generation}
  \item solve the mechanics model \eqref{eq:mechanics} coupled with the circulation model \eqref{eq:circulation}, using the time discretization step $\Delta t$; \label{it:step-mechanics}
  \item interpolate the deformation gradient $\Ftensor$ from $\mesh\mech$ to $\mesh\ep$, using interpolation between quadrature nodes (\cref{sec:interpolation-deformation-gradient,alg:rbf-svd}). \label{it:interpolation-F}
\end{enumerate}
We refer the interested reader to \cite{regazzoni2022cardiac} for more details on the methods for the electrophysiology, force generation, mechanics and circulation problems (steps \ref{it:step-ep}, \ref{it:step-force-generation} and \ref{it:step-mechanics} of the procedure above).

\begin{remark}
      In principle, we could use interpolation between quadrature nodes to interpolate $\calcium$ from the fine mesh $\mesh\ep$ to the coarse one $\mesh\mech$ (step \ref{it:interpolation-calcium} of the procedure above). However, using interpolation between \acp{DoF} is more convenient and computationally efficient, since both the ionic model and the force generation model are solved on the \acp{DoF} of the respective mesh \cite{regazzoni2022cardiac}.
\end{remark}

\subsection{Alternative interpolation methods}

The interpolation method previously described will be referred to as \rbffsvd{}. For comparison, we will also consider the following alternative interpolation schemes (each of these approaches replaces step \ref{it:interpolation-F} of the procedure described in previous section):
\begin{itemize}
      \item \nestedd{}: the meshes $\mesh\ep$ and $\mesh\mech$ are one nested into the other (i.e. the finer one is obtained from the coarser one through one or more refine-by-splitting steps). Then, the displacement field $\displacement$ is transferred from the coarse to the fine mesh through standard finite element interpolation, and the deformation gradient $\Ftensor$ is evaluated on quadrature nodes $\mathbf x\ep_j$. We refer to \cite{piersanti20223d,regazzoni2022cardiac} for additional details;
      \item \rbfd{}: the displacement field $\displacement$ is interpolated from the coarse to the fine mesh, using RBF interpolation between DoFs (\cref{sec:interpolation-dofs}); then, the deformation gradient $\Ftensor$ is evaluated on quadrature nodes $\mathbf x\ep_j$ through standard finite element interpolation. We refer to \cite{salvador2020intergrid} for further details;
      \item \rbffew{}: we perform interpolation between quadrature nodes (\cref{sec:interpolation-quadrature}) of each cartesian component of the deformation gradient $\Ftensor$; this is also known as Euclidean interpolation \cite{satheesh2022structure}.
\end{itemize}

\section{Numerical experiments}
\label{sec:numerical-experiments}

\begin{figure}
      \centering

      \includegraphics[width=0.65\textwidth]{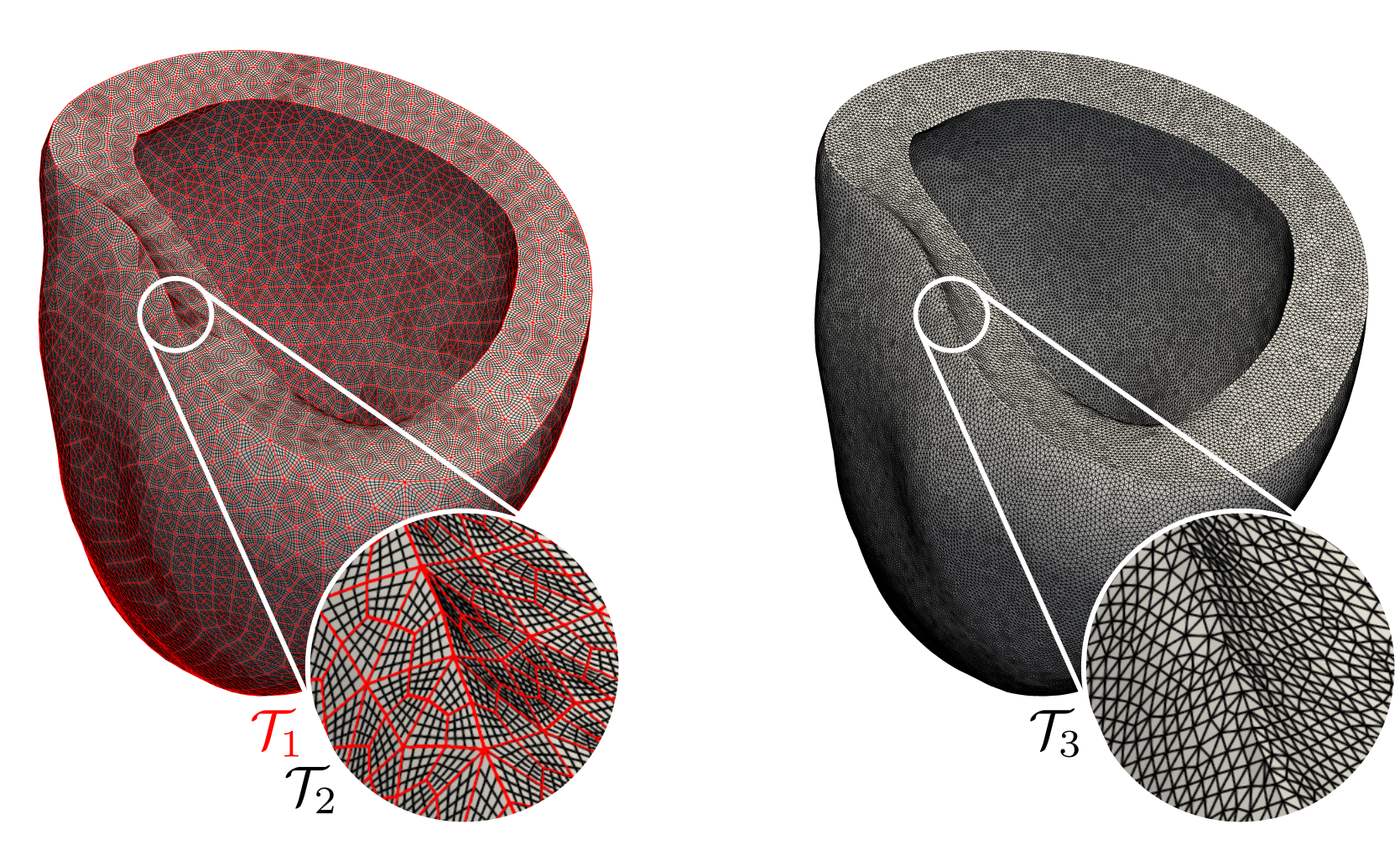} 

      \caption{Left: the \revb{meshes $\mesh_1$ (red) and} $\mesh_2$ (red). Right: the mesh $\mesh_3$.}
      \label{fig:mesh} 
\end{figure}

All simulations discussed below are performed using \lifex{}, a C++ high-performance library for cardiac applications \cite{africa2022lifexcore,africa2023lifexcfd,africa2023lifexfiber} based on the finite element core \dealii{} \cite{arndt2020dealii,arndt2022dealii9.4}. The computational domain is the left ventricle of the Zygote Heart Model \cite{zygote}, pre-processed using the techniques described in \cite{fedele2021polygonal}. We report in \cref{tab:mesh} details on the meshes (see \cref{fig:mesh}) used throughout the numerical experiments discussed below.

All simulations presented in the following sections are run on the GALILEO100 supercomputer\footnote{Technical specifications: \url{https://wiki.u-gov.it/confluence/display/SCAIUS/UG3.3\%3A+GALILEO100+UserGuide}.} from the CINECA supercomputing center (Italy). Unless otherwise specified, simulations are run in parallel using \num{192} cores.

\begin{table}
      \centering

      \begin{tabular}{c c c c c c c}
            \toprule
            & & & & \multicolumn{3}{c}{\textbf{element diameter} [\si{\milli\metre}]} \\
            \textbf{mesh} & \textbf{element type} & \textbf{\# elements} & \textbf{\# nodes} & \textbf{min.} & \textbf{avg.} & \textbf{max.} \\
            \midrule
            $\mesh_1$ & hexahedra & \num{140120} & \num{159443} & \num{1.1} & \num{2.0} & \num{3.5} \\
            $\mesh_2$ & hexahedra & \num{8967680} & \num{9141953} & \num{0.2} & \num{0.5} & \num{1.4} \\
            $\mesh_3$ & tetrahedra & \num{3838394} & \num{762797} & \num{0.4} & \num{0.8} & \num{1.2} \\
            \bottomrule
      \end{tabular}

      \caption{Type and number of elements, number of nodes, and element diameter (minimum, average and maximum) of the meshes used in the numerical experiments.}
      \label{tab:mesh}
\end{table}

\begin{table}
      \centering

      \begin{subtable}{0.45\textwidth}
            \centering

            \begin{tabular}{c S[table-format=4.3] l}
                  \toprule
                  \textbf{Parameter} & \textbf{Value} & \textbf{Unit} \\
                  \midrule
                  $\sigmaf$ & 1.68e-4 & \si{\square\metre\per\second} \\
                  $\sigmas$ & 0.769e-4 & \si{\square\metre\per\second} \\
                  $\sigman$ & 0.248e-4 & \si{\square\metre\per\second} \\
                  $T_\text{act}^\text{max}$ & 600 & \si{\kilo\pascal} \\
                  \bottomrule
            \end{tabular}

            \caption{}
            \label{tab:parameters-physical}
      \end{subtable}
      \begin{subtable}{0.45\textwidth}
            \centering

            \begin{tabular}{c c c c}
                  \toprule
                  \textbf{interpolated quantity} & $M$ & $\alpha$ \\
                  \midrule
                  $\calcium$ & 1 & 2.5 \\
                  $\displacement$ & 5 & 3.0 \\
                  $\Ftensor$ (Euclidean) & 2 & 2.0 \\
                  $\Ftensor$ (SVD) & 2 & 2.0 \\
                  \bottomrule
            \end{tabular}

            \caption{}
            \label{tab:parameters-rbf}
      \end{subtable}

      \caption{(a) Physical parameters for the simulations of \cref{sec:comparison}. Conductivities are calibrated to obtain conduction velocities along fibers, sheets and cross-fibers of \SI{0.6}{\metre\per\second}, \SI{0.4}{\metre\per\second} and \SI{0.2}{\metre\per\second} on the mesh used for electrophysiology \cite{augustin2016anatomically}. We only report parameters whose value is different from that used in \cite{regazzoni2022cardiac}. (b) Parameters controlling the adaptive RBF radius in the tests of \cref{sec:comparison,sec:flexibility}.}
\end{table}

\subsection{A comparison of interpolation methods}
\label{sec:comparison}

As a first test, we use $\mesh_1$ as a coarse mesh, for mechanics and force generation, and $\mesh_2$ as a fine mesh, for electrophysiology (see \cref{tab:mesh}). The two meshes are one nested within the other, that is the mesh $\mesh_2$ is obtained by applying a refine-by-splitting procedure to $\mesh_1$ two times (thus subdividing every element into \num{64} smaller elements, see \cref{fig:mesh}). This choice allows to compare all four interpolation methods described above (\nestedd{}, \rbfd{}, \rbffew{} and \rbffsvd{}). We consider bilinear finite elements for all the problems involved.

We perform the simulation of a heartbeat in healthy conditions, setting $T = \SI{800}{\milli\second}$, $\Delta t = \SI{1}{\milli\second}$ and $n\ep = 10$ (so that $\Delta t\ep = \SI{0.1}{\milli\second}$).
\Cref{tab:parameters-physical} lists the value of the parameters used in setting up the model. For the sake of brevity, we only report the parameters whose value is different from that used in \cite{regazzoni2022cardiac}. 

For the RBF interpolation of $\calcium$, $\displacement$ and $\Ftensor$, we manually tuned $M$ and $\alpha$ to minimize oscillations in the interpolated field, while keeping the RBF support radius as small as possible. For RBF interpolation between quadrature points, we test both $q = \num{1}$ (i.e. the interpolation points are the barycenters of each element) and $q = \num{2}$ (corresponding to \num{8} interpolation points per element). The parameters used to select the adaptive RBF radius for the different interpolation methods in this comparison are reported in \cref{tab:parameters-rbf}.

\begin{figure}
      \centering

      \begin{subfigure}{0.3\textwidth}
            \includegraphics[width=\textwidth]{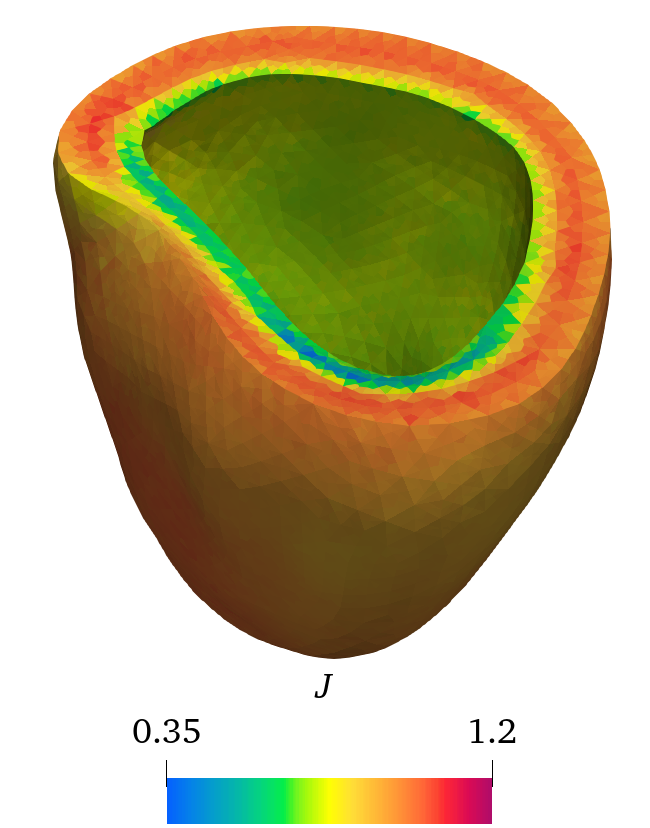}
            \caption*{coarse mesh}
      \end{subfigure}
      \begin{subfigure}{0.3\textwidth}
            \includegraphics[width=\textwidth]{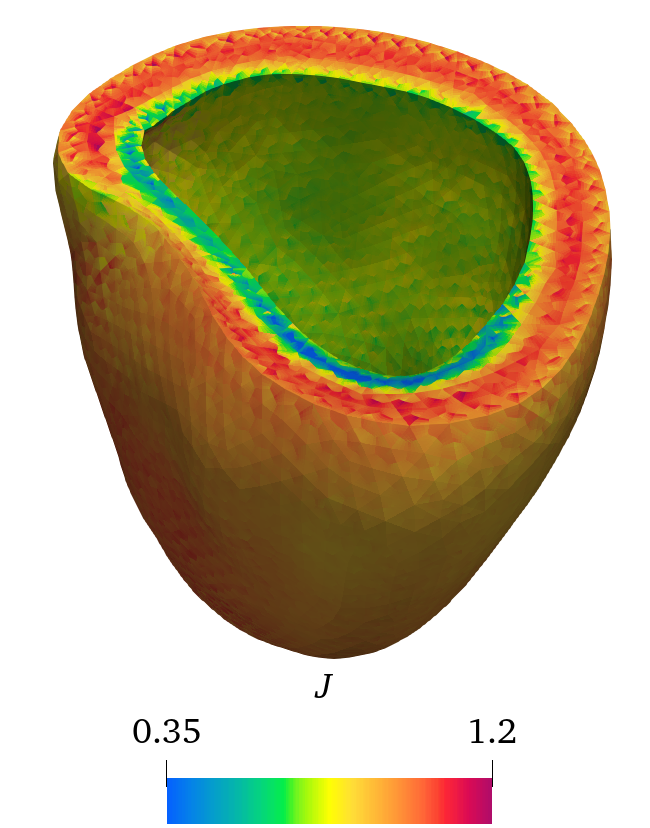}
            \caption*{\nestedd{}}
      \end{subfigure}
      \begin{subfigure}{0.3\textwidth}
            \includegraphics[width=\textwidth]{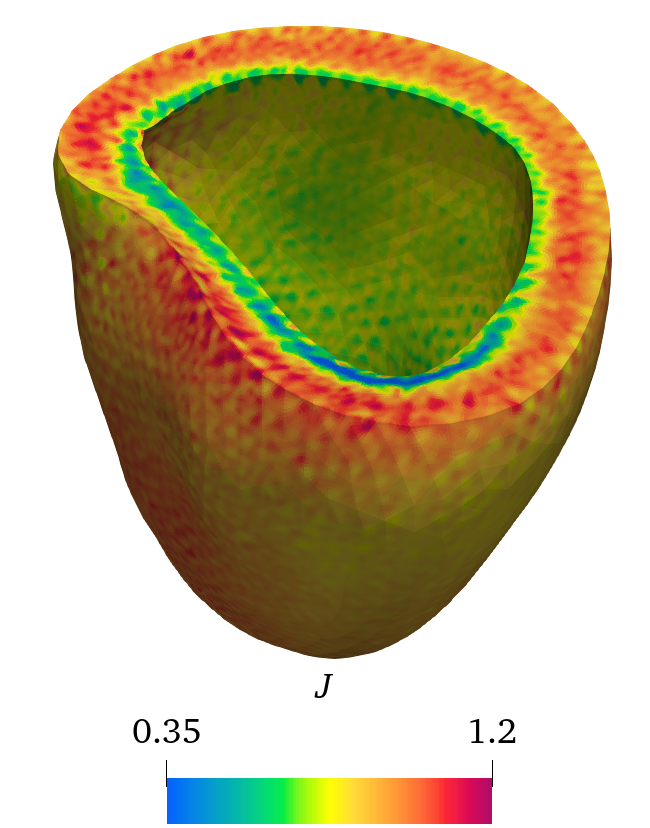}
            \caption*{\rbfd{}}
      \end{subfigure}

      \caption{For the test of \cref{sec:comparison}, determinant of the deformation gradient $\jacobian$ on the coarse mesh $\mesh_1$ (left) and on the fine mesh $\mesh_2$ using the interpolation methods \nestedd{} (center) and \rbfd{} (right). The plots are done at time $t = \SI{155}{\milli\second}$.}
      \label{fig:jacobian-3d}
\end{figure}

In this setting, the interpolation methods \nestedd{} (used in \cite{piersanti20223d,regazzoni2022cardiac}) and \rbfd{} (used in \cite{salvador2020intergrid}) fail to complete the simulation. Indeed, both yield negative values of $\jacobian$ on some of the quadrature nodes of the fine mesh, and this causes the electrophysiology solver to diverge at times $t = \SI{158}{\milli\second}$ and $t = \SI{213}{\milli\second}$, respectively. This can be explained by looking at the spatial distribution of $\jacobian$ after interpolation, as reported in \cref{fig:jacobian-3d}: although the overall distribution of $\jacobian$ is captured on the fine mesh, spurious oscillations are introduced. It is our experience that whether or not the solver crashes is extremely sensitive to the physical and numerical parameters (including the space and time discretizations). Nonetheless, the solver failure in this setting is exemplar of the unreliability of interpolation methods that transfer the displacement field and then compute its gradient.

\begin{figure}
      \centering

      \includegraphics{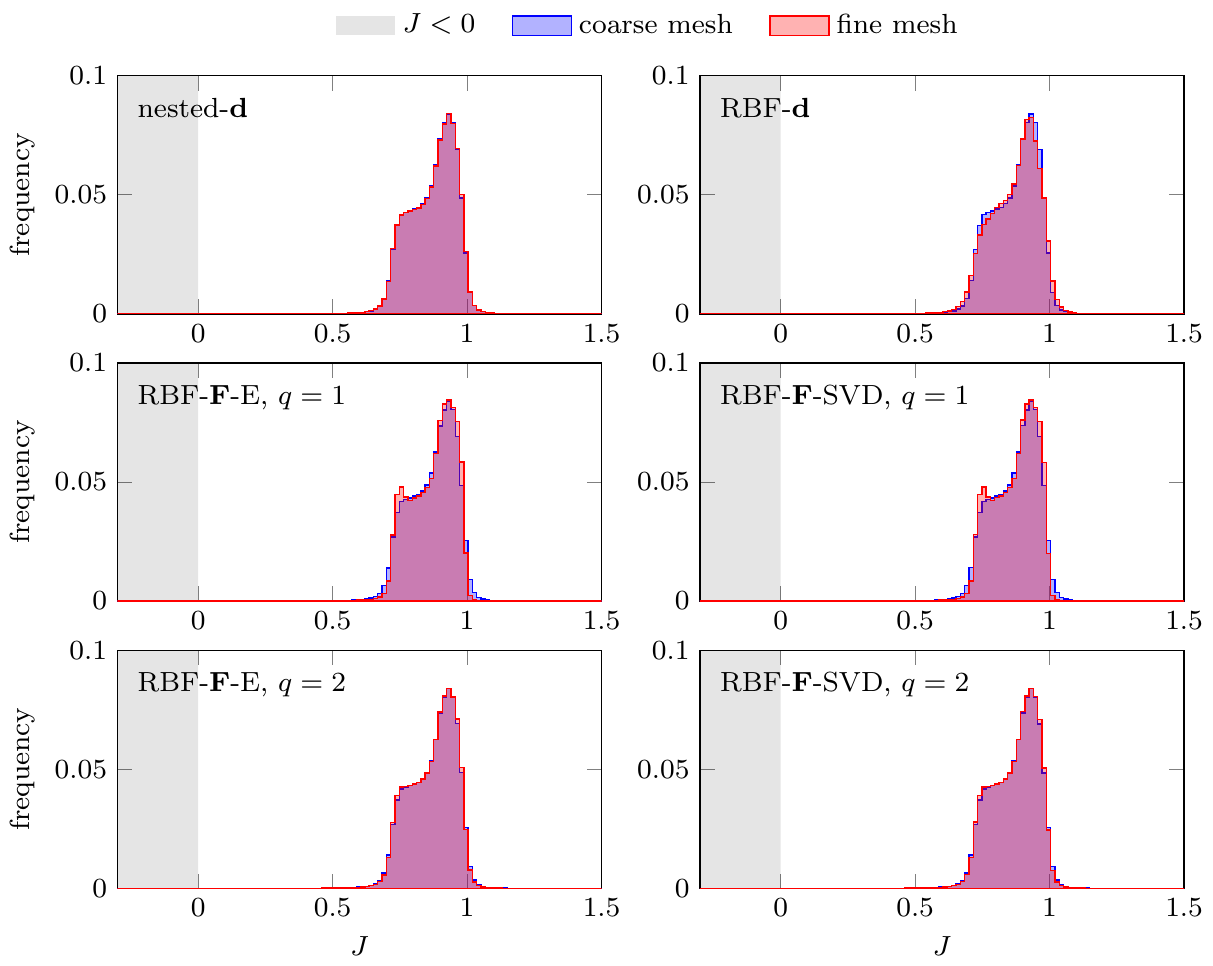}

      \caption{Histograms of the evaluation of $\jacobian$ at mesh quadrature nodes, at $t = \SI{155}{\milli\second}$, on the coarse mesh $\mesh\mech$ (blue) and on the fine mesh $\mesh\ep$ (red), with the different interpolation methods.}
      \label{fig:jacobian-histograms-linear}
\end{figure}

On the contrary, both the \rbffew{} and \rbffsvd{} interpolation methods allow the simulation to reach the final time. We attribute this difference in behavior to the fact that the interpolation methods \nestedd{} and \rbfd{} lead to values of $\jacobian$ on the fine mesh that are significantly different from those on the coarse mesh, including negative values, whereas \rbffew{} and \rbffsvd{} allow a more accurate interpolation and yield $\jacobian > 0$, ensuring the well-posedness of the discrete monodomain equation.

\begin{figure}
      \centering

      \includegraphics{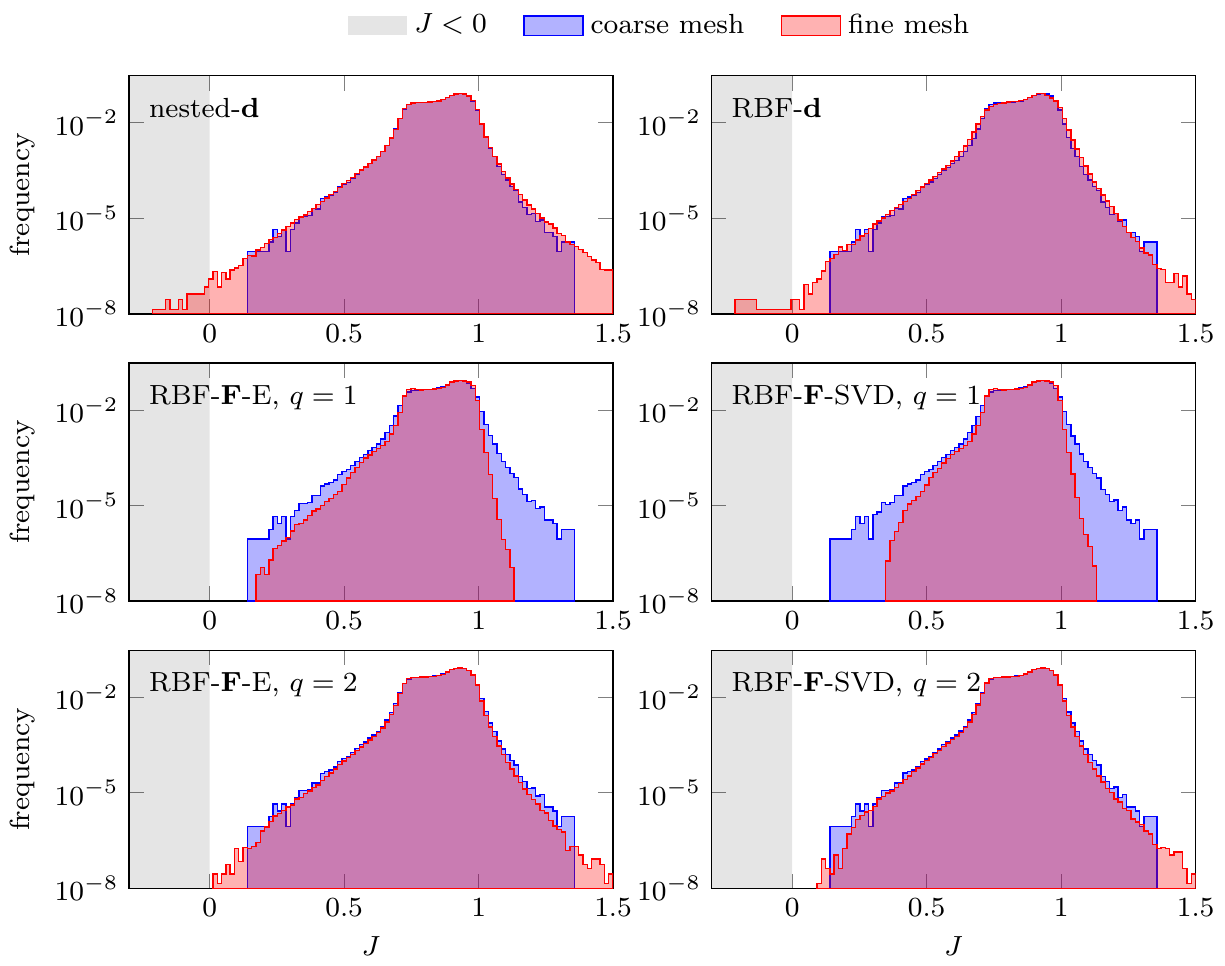}

      \caption{Histograms of the evaluation of $\jacobian$ at mesh quadrature nodes, at $t = \SI{155}{\milli\second}$, on the coarse mesh $\mesh\mech$ (blue) and on the fine mesh $\mesh\ep$ (red), with the different interpolation methods. The vertical axes are in logarithmic scale.}
      \label{fig:jacobian-histograms-logscale}
\end{figure}

To quantify this effect, we report in \cref{fig:jacobian-histograms-linear,fig:jacobian-histograms-logscale} the histograms of the values of $\jacobian$ on the quadrature nodes of $\mesh\ep$, computed with the different interpolation methods, for a representative time instant during systolic contraction ($t = \SI{155}{\milli\second}$). Although all four methods capture the overall distribution of $\jacobian$ (as seen in \cref{fig:jacobian-histograms-linear}), significant differences are present in the tails of the distributions (as highlighted by the logarithmic scale of \cref{fig:jacobian-histograms-logscale}). From the plots, we can \reva{notice} how the evaluations of $\jacobian$ with the methods based on the displacement (\nestedd{} and \rbfd{}) yield values \reva{that fall outside of the range} observed on the coarse mesh. These values are physically unrealistic (given the near-incompressibility of the solid constitutive law), and occasionally even negative (to which we attribute the failure of the solver). On the contrary, interpolation methods based on the deformation gradient (\rbffew{} and \rbffsvd{}) yield a range for $\jacobian$ that is closer to the one on the coarse mesh, in particular avoiding negative values and thus preventing the failure of the solver.

We also observe that setting $q = 1$ in the methods \rbffew{} and \rbffsvd{} has a regularizing effect, filtering out the most extreme values of $\jacobian$. On the other hand, when setting $q = 2$, thus increasing the number of interpolation points, the distribution of $\jacobian$ is more accurately recovered. This is especially true with the \rbffsvd{} scheme, whereas the \rbffew{} scheme yields a small number of points for which $\jacobian$ is close to zero. Seeing as this approach does not provide a theoretical guarantee that $\jacobian > 0$, it is possible that under different simulation settings (e.g. with a higher contractility, or in pathological conditions) this may lead to the failure of the solver. Conversely, the \rbffsvd{} method guarantees \revb{by construction} that $\jacobian > 0$ after interpolation (regardless of the choice of $q$), and as such it provides an accurate and robust tool to implement mechano-electrical feedbacks.

\begin{figure}
      \centering

      \includegraphics{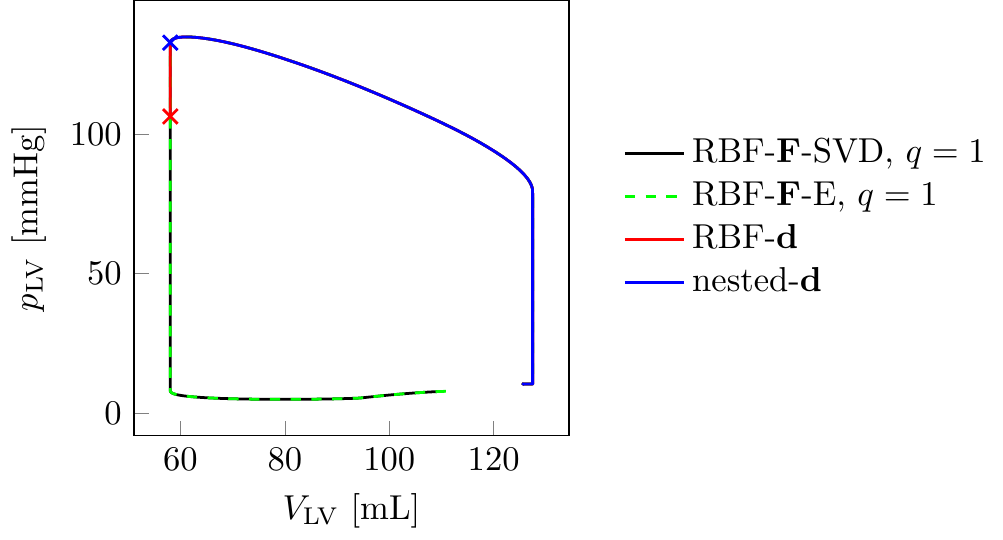}

      \caption{Pressure-volume loops of obtained with different interpolation methods for the test case of \cref{sec:comparison}. Cross markers indicate the failure of the solver. For readability, we only report the cases with $q = 1$, since those with $q = 2$ yield identical results.}
      \label{fig:pvloop-comparison}
\end{figure}

Finally, we report in \cref{fig:pvloop-comparison} the pressure-volume loops associated with all simulations. We observe that, regardless of the procedure employed, ventricular pressure and volume are essentially identical. \reva{This is in agreement} with the histograms of \cref{fig:jacobian-histograms-linear}\reva{, which show that all the methods yield comparable results for the overall distribution of $J$}. Thus, the interpolation methods based on $\Ftensor$ do not introduce perturbations with respect to the ones based on $\displacement$ and previously used in validated electromechanical simulations \cite{regazzoni2022cardiac,salvador2020intergrid,piersanti20223d}, while improving over those in terms of robustness.

\subsection{Parallel implementation and scalability}
\label{sec:implementation}

Cardiac electromechanical models are usually \reva{implemented in a parallel computing framework}, distributing the computational load across multiple processors to keep the overall time-to-solution small in spite of the large amount of unknowns. It is therefore important that the chosen intergrid interpolation method is efficient and scalable in a parallel computing setting.

\begin{figure}
      \centering

      \includegraphics{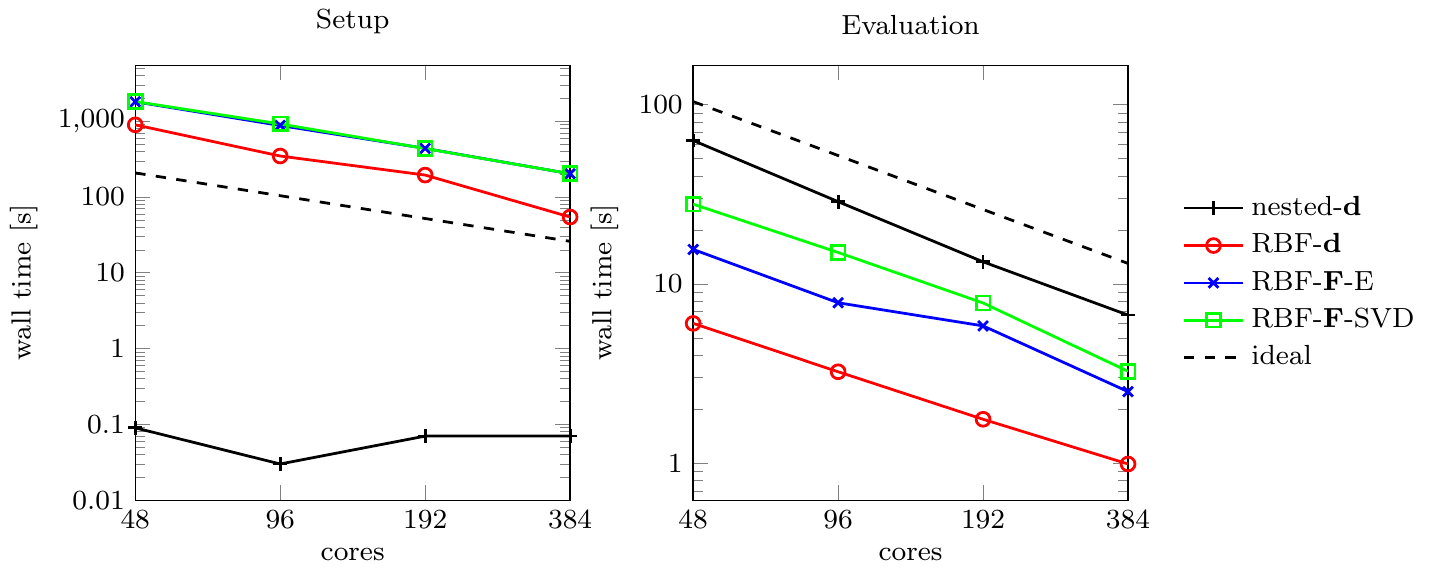}

      \caption{Strong scalability test for the coarse-to-fine interpolation. Left: wall time for the initialization of the data structures for interpolation. Right: wall time spent in evaluating the interpolant during the simulation (with a final time $T = \SI{10}{\milli\second}$).}
      \label{fig:scalability}
\end{figure}

To quantify the parallel performance of our implementation, we performed a strong scalability test for the four interpolation methods discussed before, using the same setting as in \cref{sec:comparison}. The results are reported in \cref{fig:scalability}, where we separately report the wall time spent in the initialization and the evaluation of the interpolant. We observe that, while \ac{RBF} interpolation has a significantly higher initialization cost than the nested approach (based on intergrid transfer operators implemented in \dealii{} \cite{arndt2020dealii}), the cost for initialization reduces linearly with the number of processors employed. Moreover, the initialization step consists in assembling the matrices $\interpmatrix$, $\evalmatrix$ and $\precmatrix^{-1}$. Since all of these only depend on the location of the points $\srcpoint{i}$ and $\dstpoint{i}$, this computation can be performed only once in an offline phase, and subsequently reused in multiple simulations on the same mesh, thus amortizing the initialization cost.

For all interpolation methods we observe linear scalability of the computational cost associated with the evaluation of the interpolant (\cref{fig:scalability}, right). We also observe that \ac{RBF} interpolants exhibit a slightly better performance than nested interpolation.

\begin{figure}
      \centering

      \includegraphics{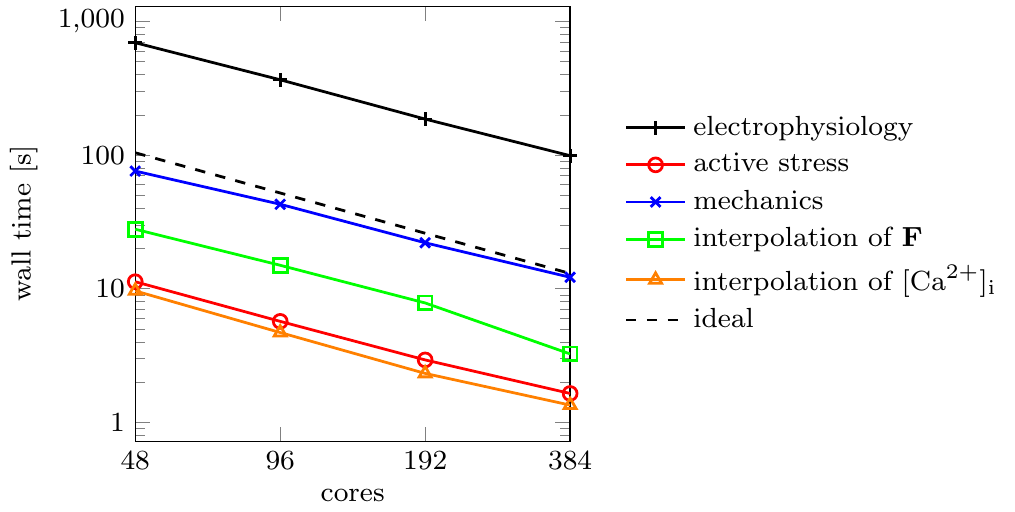}

      \caption{Comparison of the computational cost associated with intergrid interpolation and with the solution of model equations (with final time $T = \SI{10}{\milli\second}$). $\Ftensor$ is interpolated with the \rbffsvd{} method.}
      \label{fig:efficiency}
\end{figure}

Finally, \cref{fig:efficiency} reports a comparison of the computational cost associated with the evaluation of the interpolant and with the solution of the model equations themselves. The former is negligible with respect to the solution of the electrophysiology and mechanics equations (which dominate the overall computational time). We conclude that the proposed intergrid transfer method does not significantly affect the overall computational cost of the simulation.

\subsection{Flexibility of the interpolation method}
\label{sec:flexibility}

To further \revb{highlight} the flexibility of the proposed interpolation method, we consider a new test using mesh $\mesh_1$ for mechanics (as done in previous sections), and $\mesh_3$ for electrophysiology (see \cref{tab:mesh,fig:mesh}). The fine mesh is not obtained through a refine-by-splitting procedure from the coarse mesh, as done in previous examples. Instead, it is a fully independent mesh, composed of tetrahedral elements (whereas the coarse mesh $\mesh_1$ is composed of hexahedra). The electrophysiology equations are discretized using quadratic finite elements (for a total of \num{5224245} degrees of freedom). We use the \rbffsvd{} scheme for the interpolation of $\Ftensor$, setting $q = 1$ and using the same parameters as in \cref{sec:comparison} to select the RBF support radius (see \cref{tab:parameters-rbf}). We point out that the element type and the polynomial degree used for mechanics and electrophysiology are different. We use the same physical parameters as in \cite{regazzoni2022cardiac}, except for setting $T_\text{act}^\text{max} = \SI{600}{\kilo\pascal}$.

\begin{figure}
      \centering

      \begin{subfigure}[b]{0.35\textwidth}
            \centering
            \includegraphics{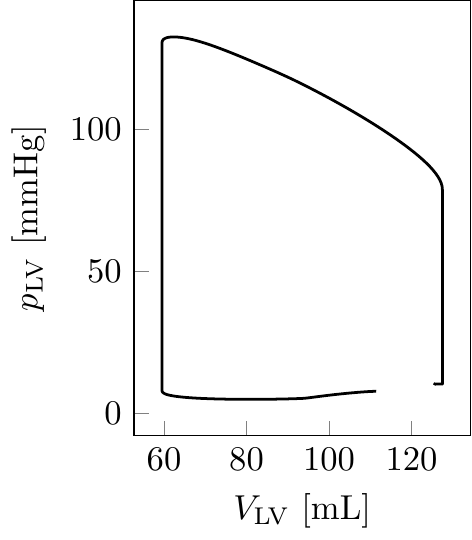}

            \caption{}
      \end{subfigure}
      \begin{subfigure}[b]{0.55\textwidth}
            \centering

            \includegraphics[width=\textwidth]{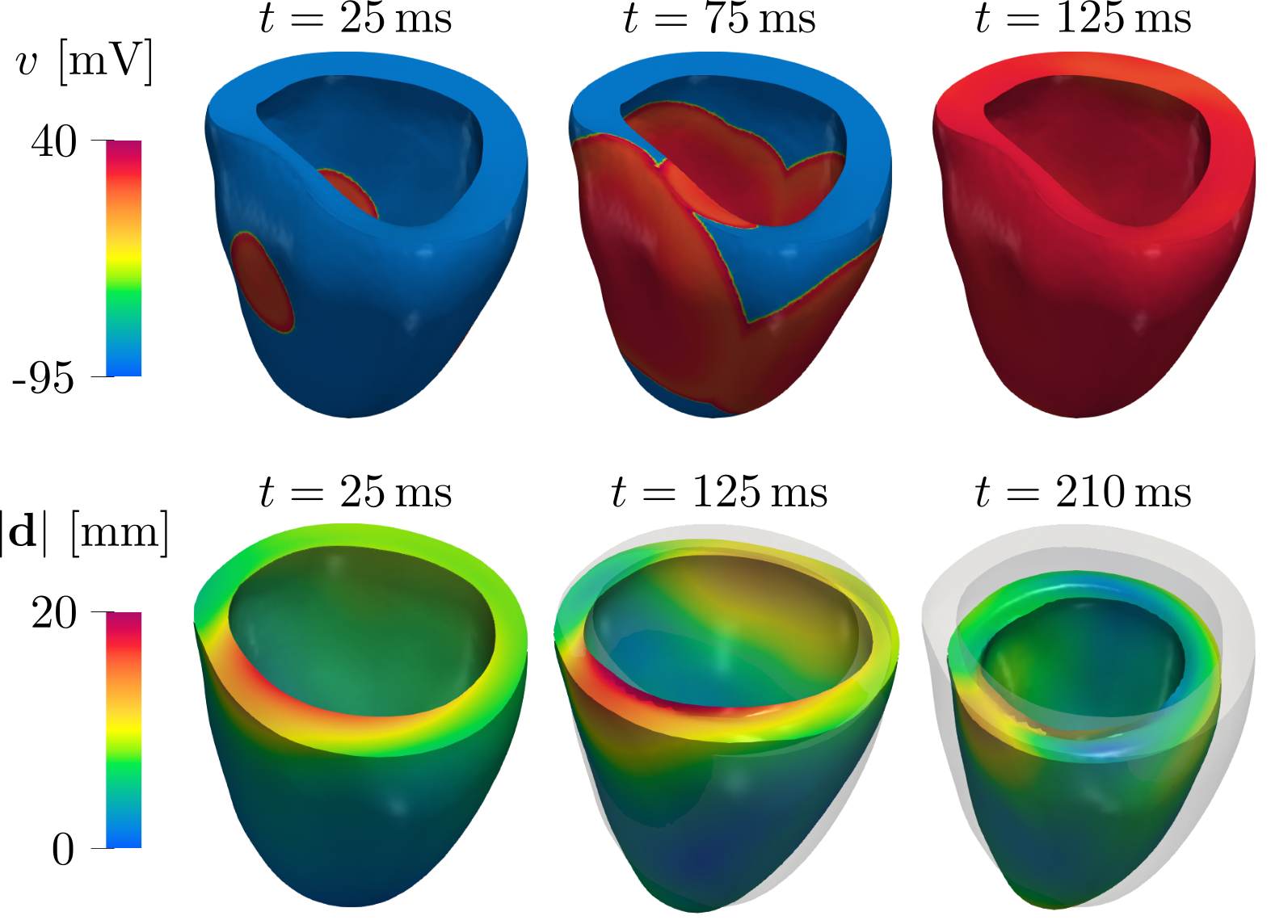}
            
            \caption{}
      \end{subfigure}

      \caption{(a) Pressure-volume loop for the simulation of \Cref{sec:flexibility}. (b) Snapshots of the transmembrane potential $\potential$ (top) and of the displacement magnitude $|\displacement|$ (bottom) for the simulation of \Cref{sec:flexibility}. The bottom plots are warped by $\displacement$, and the initial configuration is superimposed in transparency.}
      \label{fig:results-flexibility}
\end{figure}

We report the pressure-volume loop for this test case, as well as some snapshots of the solution, in \cref{fig:results-flexibility}. The results are consistent with those of \cref{sec:comparison} (as well as with those of \cite{regazzoni2022cardiac}), even if using two entirely independent (not nested) meshes for electrophysiology and mechanics.

This test shows how the proposed interpolation approach allows the transfer of solutions between radically different finite element spaces. Thus, it provides a useful tool to efficiently tackle the multiphysics nature of computational models of the heart.

\section{Conclusions}
\label{sec:conclusions}

We introduced a method to transfer the deformation gradient \reva{from a} coarse \reva{to a} fine mesh in electromechanical simulations. This is a key ingredient \revb{when using} electromechanical models \revb{that account for} mechano-electrical feedback effects (both geometrical and physiological, such as stretch-activated currents)\revb{. These models are essential to ensure accurate simulations}, especially in pathological scenarios \cite{salvador2022role}. The proposed method is based on combining rescaled, localized RBF interpolation with the SVD factorization of the deformation gradient tensor $\Ftensor$.

The proposed method is compared to alternative existing approaches in the literature, based on nested intergrid interpolation or RBF interpolation of the displacement field. Both these techniques have been previously applied to physiological simulation of cardiac electromechanics \cite{regazzoni2022cardiac, salvador2020intergrid, piersanti20223d}. The numerical experiments carried out in this work highlight the shortcomings of the existing methods, especially in terms of their robustness, by considering a numerical setting under which they lead to the failure of the electromechanical solver. We attribute this failure to the fact that the existing methods do not preserve the \reva{positivity} of $\jacobian = \det\Ftensor$, and indeed we observe points for which $\jacobian < 0$ in our numerical experiments, just before the solver failure.

Conversely, the proposed method guarantees that, after interpolation, $\jacobian > 0$, and therefore it ensures the robustness of the solver. Moreover, \reva{our results indicate} that methods based on interpolating the deformation gradient, rather than the displacement field as previously done, do not introduce artificially (and unphysically) large or small values for $\jacobian$. We believe \reva{that} this feature \reva{is} especially relevant if physiological mechano-electric feedbacks (such as stretch-activated currents) or pathological scenarios are considered. An analysis of the performance of the method under pathological conditions will be the subject of future studies.

Finally, we highlighted through a numerical experiment how the proposed method allows to easily transfer solution variables between finite elements of different degree and even of different element shape (tetrahedral or hexahedral, in our case). This result indicates that the proposed interpolation method can be a valuable tool in improving not only the robustness, but also the geometrical and parametric flexibility of electromechanical simulations, allowing to independently tailor the discretization of each core model to its specific accuracy needs.

\section*{Acknowledgements}

The authors acknowledge their membership to the GNCS - Gruppo Nazionale per il Calcolo Scientifico (National Group for Scientific Computing, Italy).
This project has been partially supported by the INdAM-GNCS Project CUP\textunderscore E55F22000270001.
LD and AQ received funding from the Italian Ministry of University and Research (MIUR) within the PRIN (Research projects of relevant national interest 2017 ``Modeling the heart across the scales: from cardiac cells to the whole organ'' Grant Registration number 2017AXL54F).
We acknowledge the CINECA award under the ISCRA initiative, for the availability of high performance computing resources and support

\clearpage
\bibliographystyle{abbrv}
\bibliography{bibliography}

\end{document}